%% file: encX.tex
\documentclass[english,12pt]{article}

\usepackage[T1]{fontenc}
\usepackage{fancyheadings}
\usepackage[utf8]{inputenc}
\usepackage{babel}
\usepackage{float}
\usepackage{mathtools}
\usepackage{amssymb}
\usepackage{bbm}
\usepackage{setspace}
\usepackage[unicode=true,pdfusetitle,
 bookmarks=true,bookmarksnumbered=false,bookmarksopen=false,
 breaklinks=false,pdfborder={0 0 0},pdfborderstyle={},backref=false,colorlinks=false]
 {hyperref}

\usepackage{tikz}
\usepackage{tikzit}
 \usepackage{tikz-cd}
 \usetikzlibrary{decorations.pathmorphing}

\usetikzlibrary{decorations,arrows}
\usetikzlibrary{decorations.markings}
\usetikzlibrary{decorations.pathmorphing}
\usetikzlibrary{matrix,arrows}

\input{sample.tikzstyles}

\everymath{\displaystyle}

\makeatletter

\begin{document}



\thispagestyle{empty}

\begin{flushright}
   {\sf ZMP-HH/23-7}\\
   {\sf Hamburger$\;$Beitr\"age$\;$zur$\;$Mathematik$\;$Nr.$\;$941}
   \\[2mm] May 2023
\end{flushright}

\vskip 2.2em

\begin{center}
{\bf \Large Algebraic structures \\[4pt] in two-dimensional conformal field theory}

\vskip 15mm

{\large \  \ J\"urgen Fuchs\,$^{\,a},~~$ Christoph Schweigert\,$^{\,b},~~$
Simon Wood\,$^{\,c,b}~~$ and $~~$ Yang Yang\,$^{d}$}

\vskip 12mm

 \it$^a$
 Teoretisk fysik, \ Karlstads Universitet\\
 Universitetsgatan 21, \ S\,--\,651\,88\, Karlstad
 \\[9pt]
 \it$^b$
 Fachbereich Mathematik, \ Universit\"at Hamburg\\
 Bereich Algebra und Zahlentheorie\\
 Bundesstra\ss e 55, \ D\,--\,20\,146\, Hamburg
 \\[9pt]
 \it$^c$
 School of Mathematics, Cardiff University\\
 Abacws, Senghennydd Road, Cardiff, CF24 4AG, UK
 \\[9pt]
 \it$^d$
Erwin Schr\"odinger International Institute for Mathematics and Physics \\
Boltzmanngasse 9, \ A\,--\,1090\, Wien

\end{center}

\vskip 3.2em

\noindent{\sc Abstract}\\[3pt]
This is an invited contribution to the 2nd edition of the Encyclopedia of 
Mathematical Physics. We review the following algebraic structures which appear
in two-dimensional conformal field theory (CFT):
 \\
The symmetries of CFTs can be
formalised as chiral algebras, vertex operator algebras or nets of observable
algebras. Their representation categories are abelian categories having additional 
structures, which are induced by properties of conformal blocks, i.e.\ of vector 
bundles over the moduli space of curves with marked points, which can be 
constructed from the symmetry structure. 
 \\
These mathematical notions pertain to the description of chiral CFTs.
In a full local CFT one deals in addition with correlators, which are specific
elements in the spaces of conformal blocks.  In fact, a full CFT is the same as
a consistent system of correlators for arbitrary conformal surfaces with any 
number and type of field insertions in the bulk as well as on boundaries and on 
topological defect lines. We present algebraic structures that allow one to 
construct such systems of correlators.

\newpage


\bigskip

\noindent
{\bf Key points}
\begin{itemize}
\item
Algebraic structures describing chiral symmetries in two-dimensional 
conformal field theories: vertex operator algebras and their representations.
\item
Structure on the representation categories of nice vertex algebras: 
braided monoidal tensor category with Grothendieck-Verdier duality.
\item
Conformal blocks as vector bundles with flat connection; 
their monodromies are encoded in terms of a modular functor that comes from a
modular fusion category.
\item
Description of consistent sets of correlators in terms of three-manifolds with 
boundary and in terms of string nets, for world sheets of arbitrary genus,
with arbitrary conformal boundary conditions and with insertions of all types of
fields, including generalised defect fields.
\end{itemize}


\section*{Introduction}

Conformal field theory in two dimensions -- CFT, for short --
is a prime example of the fruitful interplay between mathematics and physics.
On the mathematical side it combines concepts and tools from algebraic geometry, 
higher structures, modular forms, quantum topology and representation theory. 
For applications to physics, it is crucial to appreciate that there is in fact
no single notion of CFT; different applications require different theoretical setups.
For instance, conformal field theories can be defined on two-dimensional manifolds
with a metric of either Euclidean or Lorentzian signature. The qualifier
\emph{conformal} expresses the fact that the relevant geometry is the one of 
conformal manifolds, i.e.\ manifolds endowed with an equivalence class, with respect
to local rescalings, of metrics.
Working with Lorentzian signature has the advantage that structures familiar 
from local quantum field theory can be utilized \cite{frrs2,BKlr}.
On the other hand, conformal field theories defined on Euclidean manifolds have a
broader range of applications and have access to powerful tools from complex geometry.
Accordingly we focus our attention on CFTs on compact Euclidean two-manifolds.
The relation between the Euclidean and Lorentzian approaches still remains to 
be understood in full detail; however, remarkably, in both settings
similar mathematical structures arise.

A further distinction, whose importance cannot be overemphasised,
is the one between \emph{chiral} and \emph{full} conformal field theory.
Again, these are defined on two different types of two-manifolds:
   \begin{itemize}
   \item 
Together with an orientation, a Euclidean conformal structure on a two-manifold
amounts to a complex structure. \emph{Chiral} CFT is defined on a complex curve. 
Mathematically, chiral CFT consists of the theory of vertex operator algebras, their 
representations and their conformal blocks. 
   \item 
In contrast, \emph{full local CFT} is defined on conformal real two-manifolds.
Such a manifold may in particular have a non-empty boundary. In fact, one allows 
for stratified manifolds, of which manifolds with boundary are just special instances.
The one-dimensional strata of a stratified conformal two-manifold are interpreted as 
topological line defects which separate different phases of a given CFT. Specific 
types of such defects encode symmetries and dualities of CFTs \cite{ffrs5}. 
   \end{itemize}

Given a chiral CFT, a corresponding full CFT is
a \emph{consistent system of correlators}, i.e.\ a collection of 
specific elements in the spaces of conformal blocks of the chiral theory.
Both chiral and full CFTs appear in physics applications:
Chiral theories are relevant for the fractional quantum Hall effect;
in that case, the orientation has a very direct meaning as the direction of 
an external magnetic field. Full CFTs arise as world sheet theories of string theories,
in two-dimensional critical phenomena of statistical mechanics,
and in effectively one-dimensional systems in condensed matter physics.
In most of these applications one deals with oriented full CFT, i.e.\ an orientation
is chosen on the two-manifolds. But full CFTs can also defined on unoriented,
and even non-orientable, two-manifolds. Unoriented full CFTs
arise in particular as world sheet theories of type-I string theories.

   \medskip 

Best understood among all CFT models is the class of \emph{rational} CFTs, for 
which the representations of the chiral symmetries form a semisimple modular tensor
category, also known as modular fusion category; for a historical
review see \cite{fuRs13}. More recently, the focus of research has shifted to
non-rational theories, in particular to \emph{logarithmic} CFTs, for which the
conformal blocks may have logarithmic singularities and which are in particular
non-unitary; for a review of these, see e.g.\ the articles in the collection
\cite{GArR}.
By making use of powerful algebraic structures, CFT can be formulated
without recourse to a classical Lagrangian or any form of perturbation theory.
But there do exist models of (full) CFT that are based on a Lagrangian, notably
sigma-models with topological terms. Such topological terms are closely related to
bundle gerbes and other higher geometric structures, compare e.g.\ \cite{gaRe}. 

   \medskip 

This contribution is organized as follows. We first review algebraic structures
that formalize the chiral symmetries of CFTs. One conceptual framework for these
is given by conformal nets of algebras of observables \cite{BKlr}; in this framework
unitarity is deeply built in. In view of the importance of non-unitary CFTs in string
theory, statistical mechanics and as duals of supersymmetric quantum field 
theories, we restrict our attention to vertex operator algebras \cite{FRbe2} 
-- VOAs, for short -- which make no assumption on unitarity.

Afterwards we review aspects conformal blocks, which can be constructed from VOAs
and which are the building blocks of correlators
in full local CFTs. They form vector bundles with projectively
flat connections over the moduli space of complex curves. Monodromies of these
bundles can be described in terms of a modular functor. 
These data suffice to formulate the concept of a consistent set of correlators. We 
finally exhibit the construction of such consistent sets of correlators, including 
a recent approach based on a string-net formulation of modular functors.


\section*{Chiral conformal field theory}

\subsection*{Vertex operator algebras}

The basic data of a VOA are a complex vector space \(\mathfrak{V}\) and a linear map
\(Y\colon \mathfrak{V}\,{\otimes}\, \mathfrak{V} \,{\xrightarrow{~}}\, \mathfrak{V}((z))\),
called the \emph{field map}. Here \(z\) is a formal variable; in the description of
conformal blocks below it will be interpreted as a formal local coordinate on a 
complex curve. \(\mathfrak{V}((z))\) denotes \(\mathfrak{V}\)-valued power series in
$z$ and $z^{-1}$ with exponents bounded from below. The evaluation of the field map 
at a vector \(a\,{\otimes}\, b \,{\in}\, \mathfrak{V}\,{\otimes}\,\mathfrak{V}\)
is denoted by
  \begin{equation}
  Y(a,z)\, b = \sum_{n\in \mathbb{Z}}a_n b\, z^{-n-1} .
  \end{equation}
Regarding $b$ as a free variable, this is rewritten as
  \begin{equation}
  Y(a,z) = \sum_{n\in \mathbb{Z}} a_n\, z^{-n-1}
  \,\in \mathrm{End}(\mathfrak{V})[[z,z^{-1}]] \,, \quad
  a_n \,{\in}\, \mathrm{End}(\mathfrak{V}) \,.
  \label{eq:fieldser}
  \end{equation}
The axioms of a VOA then imply that \((\mathfrak{V},Y)\) is essentially a unital 
commutative associative complex algebra with a derivation and with an additional 
conformal structure. In particular, unital commutative associative complex algebras with
a derivation are examples of VOAs that do not necessarily have a conformal structure.
However, owing to the presence of formal variables, in general the commutativity and
associativity properties are more involved than in the standard commutative
algebra case. The identity element of \(\mathfrak{V}\) is called the \emph{vacuum vector}
and is frequently denoted by \(|0\rangle \); it satisfies
  \begin{equation}
  Y(|0\rangle ,z)\, a = a \qquad \text{and} \qquad
  Y(a,z)|0\rangle \in a\,{+}\,z\,\mathfrak{V}[[z]] \quad \text{for all}~ a\,{\in}\,\mathfrak{V} \,,
  \end{equation}
where \(z\,\mathfrak{V}[[z]]\) is the subspace of \(\mathfrak{V}((z))\) consisting of all
\(\mathfrak{V}\)-valued power series with only positive powers of \(z\). The conformal
structure requires the existence of a vector \(\omega\,{\in}\, V\) such that the 
coefficients $L_n$ in the expansion
  \begin{equation}
  Y(\omega,z)= \sum_{n\in \mathbb{Z}} L_n\, z^{-n-2}
  \end{equation}
satisfy the relations 
  \begin{equation}
  [L_m,L_n] = (m\,{-}\,n)\, L_{m+n} + \frac{1}{12}\, (m^3\,{-}\,m)\, \delta_{m+n,0}\, C
  \end{equation}
of the Virasoro Lie algebra. Here \(C\) is a central element; it acts on \(\mathfrak{V}\)
(and on its modules, introduced below) by multiplication with a complex number, called 
the (conformal) \emph{central charge} and denoted by \(c\). The conformal structure 
endows $\mathfrak{V}$ with an integral grading (usually required to be non-negative):
  \begin{equation}
  \mathfrak{V} = \bigoplus_{n\in \mathbb{Z}} \mathfrak{V}_n \qquad \text{with} \qquad
  \mathfrak{V}_n = \{v\,{\in}\, \mathfrak{V}\,\vert\, L_0 v \,{=}\, n v\} \,.
  \end{equation}
For a homogeneous vector \(b\,{\in}\,\mathfrak{V}_h\) it is common in the physics 
literature to shift the indexing in the series expansion \eqref{eq:fieldser} in such a way that
  \begin{equation}
  Y(b,z) = \sum_{n\in \mathbb{Z}} b_{(n)}\, z^{-n-h} .
  \end{equation}
The so obtained coefficients \(b_{(n)}\) are degree \(-n\) maps, that is, 
\(b_{(n)} \mathfrak{V}_{m} \,{\subset}\, \mathfrak{V}_{m-n}\).


\medskip

VOAs are formulated with the help of formal variables, respectively formal local 
coordinates. An algebraic structure that captures the symmetries of a chiral CFT on a
smooth algebraic curve in a coordinate-independent way is the one of a
\emph{chiral algebra} \cite{BEdr}. This structure is formulated with the
help of the notion of a $D$-module. Chiral algebras on the affine line 
$\mathbb {A}^{1}$ that are equivariant with respect to translations are in
bijection to vertex operator algebras. The study of chiral algebras has triggered
further developments, in particular it has led to the notion of a factorisation 
algebra. For detailed information we refer to \cite{BEdr,FRbe2}.

In the physics literature, the term `chiral algebra' is often used in a 
different way, as a generic term for any type of chiral symmetry structure,
and sometimes also as a synonym for VOAs.


\subsection*{Representation categories}

As a generalisation of associative algebras, VOAs naturally admit modules. A module over
a VOA \((\mathfrak{V},Y)\) consists of a complex vector space \(M\) and a linear map
\(Y^M \colon V\,{\otimes}\, M \,{\xrightarrow{~}}\, M((z))\) such
that \(Y^M\) represents the field map \(Y\) (the multiplication on
\(\mathfrak{V}\)) on \(M\). The classification of \(V\)-modules is in general an
extremely hard problem. But for interesting classes of VOAs it can be made tractable 
by restricting the kind of modules to be considered to those visible to the Zhu algebra 
formalism; see \cite{gabe15} for more information.

Continuing the analogy to commutative algebras, it is natural to study 
\(\mathfrak{V}\)-multi\-lin\-ear maps, in particular bilinear ones.
Such maps are a crucial ingredient in the construction of 
conformal blocks and are thus of central importance to chiral CFT, see below.
A \(\mathfrak{V}\)-bilinear map from a pair of
\(\mathfrak{V}\)-modules $M_1$ and $M_2$ to a third module \(M_3\) is a linear map
  \begin{equation}
  \mathcal{Y}\colon\quad M_1\,{\otimes}\, M_2 \xrightarrow{~~} M_3\{z\}[\log(z)]
  \end{equation}
that is compatible with the action of \(\mathfrak{V}\) on each of the three modules.
Here \(\log(z)\) is a formal variable satisfying \(\partial_z\log(z)\,{=}\,z^{-1}\)
and \(M_3\{z\}[\log z]\) denotes polynomials in \(\log(z)\) whose
coefficients are \(M_3\)-valued power series in $z$ for which the exponents can
be arbitrary complex numbers.

The existence of bilinear maps allows one to define a \(\mathfrak{V}\)-tensor
product, also called a \emph{fusion product} (see \cite{crkM} for a summary and 
\cite{HLZ1,HLZ2} for an exhaustive discussion). This product is characterised by
a direct generalisation of the universal property for the tensor product over a
ring: for any pair of modules $M_1$ and $M_2$ the fusion product is a
module \(M_1\,{\otimes_{\mathfrak{V}}}\, M_2\) together with an intertwining
operator \(\mathcal{Y}^{M_1,M_2}\colon M_1\,{\otimes}\, M_2 \,{\xrightarrow{~}}\,
M_1\,{\otimes_{\mathfrak{V}}}\,M_2\{z\}[\log(z)]\) such that for any module \(X\) and
intertwining operator \(\mathcal{I}\colon M_1\,{\otimes}\, M_2 \,{\xrightarrow{~}}\,
X\{z\}[\log(z)]\) there exists a unique module morphism 
\(\phi\colon M_1\,{\otimes_{\mathfrak{V}}M_2}\, \,{\xrightarrow{~}}\, X\) such that the
diagram
  \begin{equation}
  \begin{tikzcd}[column sep=1.7cm, row sep=0.9cm]
    M_1\,{\otimes}\, M_2 \ar[r,"\mathcal{Y}^{M_1,M_2}_{}"] \ar{rd}[swap]{\mathcal{I}}
    & M_1\otimes_{\mathfrak{V}} M_2\{z\}[\log(z)]
    \ar[d,"\exists! \phi",dashed]
  \\
    & X\{z\}[\log(z)]
  \end{tikzcd}
  \label{eq:univprop4otimes}
  \end{equation}
commutes. Being defined through a universal property, fusion products are unique if 
they exist. As a consequence, well chosen
categories of modules (in particular, a morphism $\phi$ making \eqref{eq:univprop4otimes}
commutative must exist within the category for every pair of objects) are furnished with 
the structure of a monoidal category, and even of a balanced braided monoidal category.
This structure can be characterised as follows:
  \begin{itemize}
\item
The tensor unit is given by the VOA \(\mathfrak{V}\). The left unit constraint $\ell$
identifies the intertwining operator \(\mathcal{Y}^{\mathfrak{V},M}\) with 
the action \(Y^M\) of \(\mathfrak{V}\) on \(M\), that is,
  \begin{equation}
  \ell_M(\mathcal{Y}^{\mathfrak{V},M}(a,z)\,m) = Y^M(a,z)m \qquad \text{for}~~
  a\,{\in}\,\mathfrak{V}\,,\ m\,{\in}\, M \,.
  \end{equation}
\item 
The associativity constraint $A$ corresponds to identifying compositions of
intertwining operators: 
  \begin{multline}
  A_{M_1,M_2,M_3}(\mathcal{Y}^{M_1,M_2\otimes_{\mathfrak{V}}M_3}(m_1,z_1)
  \mathcal{Y}^{M_2,M_3}(m_2,z_2)m_3)
  \\
  = \mathcal{Y}^{M_1\otimes_{\mathfrak{V}}M_2,M_3}(\mathcal{Y}^{M_1,M_2}(m_1,z_1{-}z_2)m_2,z_2)m_3 \,.
  \end{multline}
Geometrically, thinking of the variables $z_i$ as complex coordinates, the left 
hand side is expanded in a domain for which the insertion point $z_2$ of \(M_2\)
is close to that of \(M_3\), i.e.\ \(0\), while the right hand side is in a domain
where \(M_1\) (inserted at \(z_1\)) is close to \(M_2\). The existence of
associator isomorphisms satisfying pentagon equations is the main
obstruction for a chosen category of \(\mathfrak{V}\)-modules to be monoidal.
(The map \(A_{M_1,M_2,M_3}\) depends on the insertion points $z_1$ and $z_2$. The 
actual categorical associator is a certain limit; for details see
\cite[Sect.\,3.3]{crkM}.  For notational simplicity we suppress this issue.)
\item 
The braiding isomorphisms $c$ correspond to identifying the intertwining operator
\(\mathcal{Y}^{M_1,M_2}\) at \(z\) with the intertwining operator
\(\mathcal{Y}^{M_2,M_1}\) transported to \(-z\), that is,
  \begin{equation}
  c_{M_2,M_1}^{}(\mathcal{Y}^{M_2,M_1}(m_2,z)\,m_1)
  = \mathrm e^{zL_{-1}}\mathcal{Y}^{M_1,M_2}(m_1,-z)\, m_2 \,.
  \end{equation}
The geometric intuition for these isomorphisms is the exchange of the location
of the modules \(M_1\) and \(M_2\): 
Multiplying \(z\) by \(-1\) is a rotation of \(z\) (where \(M_1\) is located) by the 
angle \(\pi\) around \(0\) (where \(M_2\) is located), while the \(L_{-1}\)-exponential 
is a translation of both points by \(-z\), so that altogether the locations of the two
points are switched. (A choice of logarithm needs to be made for \(-1\) 
to distinguish between clockwise and counter-clockwise rotations.)
\item 
There is a twist isomorphism $\theta$ given by \(\theta_M \,{=}\,\exp(2\pi i L_0)|_M\),
which is balanced with respect to the braiding, that is,
  \begin{equation}
  \theta_{M_1\otimes_{\mathfrak{V}}M_2} = c_{M_2,M_1}^{} \,{\circ}\, c_{M_1,M_2}^{}
  \,{\circ}\, (\theta_{M_1}\,{\otimes_{\mathfrak{V}}}\, \theta_{M_2}) \,.
  \end{equation}
The twist thus defines a balanced structure on the category. In particular,
if the category is in addition rigid, then the twist defines a ribbon structure.
 \end{itemize}

There are a number of additional conditions that \(\mathfrak{V}\) or its
chosen of category of modules may satisfy, leading to the existence of
additional structure. For example, if the category is closed under taking
gradewise duals (also known as contragredient modules), then it is a ribbon
Grothendieck Verdier category (a type of monoidal category with a notion of duality 
that can be weaker than rigidity). If \(\mathfrak{V}\) satisfies a technical condition
called \(C_2\)-cofiniteness, then there is a natural choice of category (called 
admissible modules) that is finite, i.e.\ it has finitely many simple objects, every
object has finite length, and all morphism spaces are finite-dimensional. 
Finally, if \(\mathfrak{V}\) is \(C_2\)-cofinite and if the category 
$\mathrm{Rep}(\mathfrak{V})$ of admissible modules is semisimple, then it is a modular
fusion category \cite{huan21}. A VOA having this property is called a \emph{rational} 
VOA. Similar results yielding modular fusion categories can also be obtained in the 
framework of nets of observables \cite{kalM}.

A striking feature of rational VOAs is the behaviour of the characters of their
modules under modular transformations. 
The \emph{character} of a \(\mathfrak{V}\)-module \(M\) is the power series
  \begin{equation}
  \chi_M^{}(\tau) = \mathrm{tr}_M^{} q^{L_0-c/24} \quad \text{with} \quad
  q = \exp(2\pi i \tau)\,,\ \tau \,{\in}\, \mathbb{H} \,,
  \label{eq:def:char}
  \end{equation}
defined on the complex upper half plane \(\mathbb{H}\). For rational 
$\mathfrak V$, the prescriptions
  \begin{equation}
  S(\chi_M(\tau)^{}) = \chi_M^{}({-}\frac{1}{\tau}) \quad \text{and} \quad
  T(\chi_M(\tau)^{}) = \chi_M^{}(\tau{+}1) \,,
  \label{eq:SandT4characters}
  \end{equation}
called the modular \(S\)- and \(T\)-transformations, respectively, even give rise 
(after, if needed, refining the information in the
characters by considering dependence on further variables) to
an action of the group \(\mathsf{SL}(2,\mathbb{Z}) \,{\cong}\, \langle S,T \,\vert\,
S^4 \,{=}\, 1, (ST)^3 \,{=}\, S^2\rangle\), which is the modular group of the 
one-punctured torus, on the characters. More specifically, after
picking a finite set $\mathrm I \,{=}\, \{ M_i \}_i^{}$ of 
representatives for the isomorphism classes of simple modules such that 
$M_0 \,{=}\, \mathfrak{V}$, the modular transformations 
\eqref{eq:SandT4characters} are realized by numerical square matrices according to
  \begin{equation}
  S(\chi_{M_i}^{}) = \sum_{j\in \mathrm I} S_{i,j}\, \chi_{M_j}^{} \quad \text{and} \quad
  T(\chi_{M_i}^{}) = T_{i,i}\,\chi_{M_i}^{}
  \label{eq:SmatTmat}
  \end{equation}
for $i \,{\in}\, \mathrm I$.
On the other hand, as the category $\mathrm{Rep}(\mathfrak{V})$ for a rational VOA
$\mathfrak{V}$ is a modular fusion category, via the trace of a double braiding 
and via the eigenvalues of the ribbon twist, one can define categorical \(S\)- and 
\(T\)-matrices which again generate the modular group $\mathsf{SL}(2,\mathbb{Z})$. 
After a canonical rescaling, these are equal to the modular transformation matrices of
the characters defined in \eqref{eq:SmatTmat}, and as a consequence the tensor product
in $\mathrm{Rep}(\mathfrak{V})$ can be expressed through the Verlinde formula, which
depends only on the entries of the modular \(S\)-matrix \cite{huan21}. Thus in the 
rational case the decomposition rules for fusion products can be computed from the 
modular behaviour of characters. Determining the modular transformations of characters
is significantly more tractable than computing fusion products directly.

Characters are particular examples of torus chiral one-point correlators (see the 
section on conformal blocks). Thus their modular behavior is a consequence of the 
geometry of the curves they are defined on. Accordingly it is believed that 
the finiteness conditions that make up rationality are not necessary conditions;
they are, however, the only fully understood case. Still, various proposals for 
Verlinde-like formulas beyond rationality have been made, see e.g.\ \cite{riWo3}
and the literature cited there.

The information encoded in the modular fusion category $\mathrm{Rep}(\mathfrak{V})$
is often referred to as the \emph{chiral data}, or also Moore-Seiberg data, of a CFT
\cite{mose3}. The \emph{modular data} of the CFT are given by the subset of
chiral data consisting of the matrices $S$ and $T$ and the central charge $c$.


\subsection*{Conformal blocks}

\subsubsection*{Ward identities and conformal blocks}

A common theme in quantum field theory is the quest for establishing a correspondence
between states and fields. 
Furthermore, one would like to formulate the theory on a whole class of suitable
manifolds. In the present context, this amounts to the idea to associate to every 
vector $v \,{\in}\, {\mathcal H}_\lambda$ in a VOA-module ${\mathcal H}_\lambda$ 
a suitable `field' $\Phi_\lambda(v,p)$ depending on a point $p$ on complex curve $X$, 
and to assign, given any complex curve $X$ of genus $g$ and any choice of $m$ pairwise
distinct points $p_i$ on $X$, to an $m$-fold product of such fields a `correlator' 
$\langle \Phi_{\lambda_1}(v_1,p_1) \cdots \Phi_{\lambda_m}(v_m,p_m)\rangle$.
As will be seen below, formalizing this physics idea naturally leads to 
quantities which do not enjoy all properties that one would normally require for
the correlators of a local quantum field theory. Accordingly we have put 
the terms `field' and `correlator' in quotation marks.

The global chiral symmetries of the CFT that are encoded in a VOA $\mathfrak V$ imply 
linear relations that the `correlators' must satisfy. These relations are called the 
chiral \emph{Ward identities}; their solutions form vector spaces 
$B_{\vec\lambda}(X,\vec p)$, called the spaces of 
\emph{conformal blocks}. To write the Ward identities in a compact form, we regard
a conformal block as a linear form on the vector space tensor product 
  \begin{equation}
  \vec{\mathcal H}_{\vec\lambda} := {\mathcal H}_{\lambda_1}{\otimes_{\mathbb C}^{}}
  \cdots {\otimes_{\mathbb C}^{}}\, {\mathcal H}_{\lambda_m}
  \end{equation}
of $\mathfrak V$-modules that depends on the complex structure of the curve $X$ and
on the positions $\vec p \,{=}\, (p_1,p_2,...\,,p_m)$ of the field insertions,
that is, as a particular vector
  \begin{equation}
  \beta_{X,\vec p;\vec\lambda} \,\in\, \vec{\mathcal H}_{\vec\lambda}^*
  \end{equation}
in the (algebraic) dual of the tensor product. Then the `correlator' is the scalar
 
  \begin{equation}
  \langle \Phi_{\lambda_1}(v_1,p_1) \cdots \Phi_{\lambda_m}(v_m,p_m)\rangle  
  = \beta_{X,\vec p;\vec\lambda}(v_1\,{\otimes}\cdots{\otimes}\, v_m) \,.
  \end{equation}


\subsubsection*{WZW conformal blocks as invariants}

A crucial step in the construction of conformal blocks is to obtain for a given
curve $X$ a `global' variant of the VOA. For a general vertex operator algebra, this
has been analysed in terms of bundles of VOAs in Chapter 6 of \cite{FRbe2}, and in 
terms of `chiral Lie algebras' in \cite{naTs,dagT2}. For the present
purposes, we content ourselves to convey the idea by spelling it out for the arguably
simplest class of CFTs, namely the one of chiral \emph{WZW models}, as expounded in
\cite{beau8}. For these the VOA $\mathfrak V$ is generated by currents $J^a(z)$ whose 
zero modes $J^a_0$ span a finite-dimensional complex simple Lie algebra 
$\bar{\mathfrak g}$, and a $\mathfrak V$-module ${\mathcal H}_\lambda$ carries an
action of the untwisted affine Lie algebra ${\mathfrak g} \,{=}\, 
\bar{\mathfrak g}^{\scriptscriptstyle(1)}$ associated with $\bar{\mathfrak g}$.

In this situation we can combine the Lie algebra $\bar{\mathfrak g}$ with the 
commutative associative algebra ${\mathcal F}_{\!X,\vec p}$ of holomorphic 
functions on $X{\setminus}\,\vec p$ that at each of the points $p_i$ have at most
a finite order pole, which yields a Lie algebra 
$\bar{\mathfrak g}\,{\otimes}\,{\mathcal F}_{\!X,\vec p}$ of 
$\bar{\mathfrak g}$-valued functions. This \emph{global Lie algebra} captures the 
global aspects of the VOA. (For analogous global Lie algebras for general VOAs 
see \cite[Ch.\,19.4]{FRbe2} and \cite[Sect.\,3]{dagT2}.) An action of 
$\bar{\mathfrak g}\,{\otimes}\,{\mathcal F}_{\!X,\vec p}$ on the tensor product 
$\vec{\mathcal H}_{\vec\lambda}$ is obtained as follows: For each $i$ choose a
local holomorphic coordinate $\xi_i$ around the insertion point $p_i$ and expand 
$f\,{\in}\, {\mathcal F}_{\!X,\vec p}$ in a Laurent series
  \begin{equation} 
  f^{(i)}(\xi_i) = \sum_{n\gg -\infty} a_n^{(i)}\,\xi_i^n .
  \end{equation}
Then associate to the element $J^a_0\,{\otimes}\, f$ of
$\bar{\mathfrak g}\,{\otimes}\,{\mathcal F}_{\!X,\vec p}$ the element
  $
  \sum_{n\gg-\infty} a_n^{(i)} J^a_n
  $
of the affine Lie algebra ${\mathfrak g}^{\scriptscriptstyle(1)}$, acting on the
$\mathfrak V$-module ${\mathcal H}_{\lambda_i}$. The action of
$J^a_0\,{\otimes}\, f$ on $\vec{\mathcal H}_{\vec\lambda}$ is then defined as the sum
  \begin{equation}
  \sum_{i=1}^m \mathbf1 \otimes \cdots\otimes \Big( \!\sum_{n\gg-\infty}\!\!
  a_n^{(i)} J^a_n \Big) \otimes \cdots\otimes \mathbf1 \,, 
  \end{equation}
where in the $i$th summand the non-trivial tensor factor is located at the $i$th 
position.
The space of conformal blocks can now be defined as the vector space of invariants
for the induced action on the dual space $\vec{\mathcal H}_{\vec\lambda}^*$:
  \begin{equation} 
  B_{\vec\lambda}(X,\vec p) := {(\vec{\mathcal H}_{\vec\lambda}^*)}_0 \,. 
  \label{eq:block}
  \end{equation}
The action of $\bar{\mathfrak g}\,{\otimes}\,{\mathcal F}_{\!X,\vec p}$ depends 
on the choice of the local coordinates $\xi_i$, and thus the space 
$B_{\vec\lambda}(X,\vec p)$ of conformal blocks depends on that choice as 
well. However, for WZW models, and, more generally, for $C_2$-cofinite rational VOAs,
the Virasoro algebra, obtained from the conformal structure of the
VOA, provides a natural action of the group of changes of local coordinates on the 
VOA-modules, so that the conformal blocks transform covariantly under such choices.


\subsubsection*{Bundles of conformal blocks}

Next we describe what happens if the shape of the complex curve $X$ and the location 
of the insertion points $\vec p$, which have so far been kept fixed, are varied. These
data specify a point in the moduli space ${\mathcal M}_{g,m}$ of curves of genus $g$
with $m$ distinct marked points,
and together with choices of local coordinates, a point in a larger moduli space 
$\widetilde{\mathcal M}_{g,m}$. The vector spaces $B_{\vec\lambda}(X,\vec p)$ 
define a quasi-coherent sheaf on $\widetilde{\mathcal M}_{g,m}$. If the VOA is
rational and $C_2$-cofinite, then \cite[Sect.\,8.7]{dagT2} the sheaf is coherent 
and descends to the moduli space ${\mathcal M}_{g,m}$, so that the spaces 
$B_{\vec\lambda}(X,\vec p)$ fit together into the total space of a vector bundle 
$\mathcal B_{g,\vec\lambda}$ over the moduli space ${\mathcal M}_{g,m}$;
hereby in particular the dependence on choices of local coordinates is removed.
Moreover, the vector bundle $\mathcal B_{g,\vec\lambda}$ is equipped with the
additional structure of a projectively flat connection with regular singularities,
which is called the \emph{Knizhnik-Zamolodchikov connection}. This can be traced back
to the existence, for any VOA-module $\mathcal H_\lambda$, of a flat connection 
$\nabla\,{=}\, \mathrm d\,{+}\,L_{-1}\,{\otimes}\, \mathrm d z$ on a vector bundle 
$\widetilde{\mathcal H}_\lambda$ over ${\mathcal M}_{g,m}$ whose fibres are given
by $\mathcal H_\lambda$. The existence of the Knizhnik-Zamolodchikov connection 
amounts to a projective action of the fundamental group of ${\mathcal M}_{g,m}$, i.e.\
of the \emph{mapping class group} $\mathrm{Map}_{g,m}$ of the $m$-punctured curve
$X$, on each fibre, i.e.\ on the vector space $B_{\vec\lambda}(X,\vec p)$.

The term `conformal block' is frequently also used for the sheaves of local
horizontal sections in the vector bundles $\mathcal B_{g,\vec\lambda}$. 
As the bundles $\mathcal B_{g,\vec\lambda}$ are generically non-trivial, 
these sections are multivalued functions of the insertion points. Thus, as already
pointed out, they are not ordinary functions, in contrast to the correlators -- 
also called \emph{correlation functions} -- of a local quantum field theory.
Via the Knizhnik-Zamolodchikov connection, the horizontality of the sections
translates to a first order differential equation that the conformal blocks must
satisfy, called the \emph{Knizhnik-Zamolodchikov equation}.

There is also a kind of converse of the above derivation of the conformal block
spaces: From a system of conformal blocks one can construct a VOA by restricting 
to insertions of vectors that only produce a meromorphic dependence in correlation 
functions on the sphere. The axioms of a VOA and its modules as well as further 
structures can then be extracted from the desired properties of correlations 
functions. In fact, as demonstrated in \cite{gago}, to this end it is even sufficient 
to know the values of the linear forms \eqref{eq:block} only on certain
finite-dimensional 
vector spaces, which are then recognized as particular subspaces of the VOA.


\subsubsection*{Factorisation}

In the considerations above, the curve $X$ has implicitly been assumed to
be smooth. However, many arguments still go through if $X$ is instead allowed to 
be \emph{stable} and thus can
possess a mild form of singularities, known as ordinary double points. Such a
double point $p$ on $X$ can be `blown up', which results in a smooth curve
$X'$ with a projection onto $X$ under which $p$ has two pre-images $p'_\pm$.
A further crucial aspect of conformal blocks is \emph{factorisation}, which
describes their behaviour in such a situation. At present, this is thoroughly
understood only for $C_2$-cofinite rational VOAs $\mathfrak V$. In this case the 
conformal blocks even form vector bundles over the moduli spaces 
$\overline{\mathcal M}_{g,m}$ of stable pointed curves, the vector spaces of 
conformal blocks are finite-dimensional, and there exist canonical isomorphisms 
\cite[Thm.\,8.4.1]{dagT2}
  \begin{equation}
  g_{X,X'}^{}\colon \quad 
  B_{\vec\lambda}(X,\vec p) \,\xrightarrow{~\cong~}\, \bigoplus_{\mu\in I}
  B_{\vec\lambda\cup\{\mu,\mu^\vee\}}(X',\vec p \cup\! \{p'_+,p'_-\})
  \label{eq:fact}
  \end{equation}
between the spaces of conformal blocks on $X$ and $X'$, where the (finite) summation
is over the isomorphism classes of simple objects of $\mathrm{Rep}(\mathfrak V)$.
This structure tightly links the system of vector bundles 
$\mathcal B_{g,\vec\lambda}$ over the moduli spaces $\overline{\mathcal M}_{g,m}$ 
for all possible values of the genus $g$ and the number $m$ of insertion points.
 
The vector spaces $B_{\vec\lambda}(X,\vec p)$, as well as the sheaves of
horizontal sections of the bundles $\mathcal B_{g,\vec\lambda}$, depend functorially
on the VOA-modules at the insertion points, i.e.\ on the objects $\lambda_i$ of the
category $\mathrm{Rep}(\mathfrak V)$. From a categorical perspective, it is natural
to expect that the factorisation isomorphisms \eqref{eq:fact} can be generalised to
classes of non-rational VOAs, with the direct sum replaced by the categorical
notion of a \emph{coend} (similarly as in formula \eqref{eq:coendBlC} below).

An important consequence of factorisation is the \emph{Verlinde formula} for rational
VOAs, which expresses the rank of the bundles $\mathcal B_{g,\vec\lambda}$ for any
value of $g$ and $m$ through the entries of the matrix $S$ that according to
\eqref{eq:SandT4characters} describes the effect
of the modular transformation $\tau \,{\mapsto}\, {-}1/\tau$ on the characters of the 
VOA-modules $\mathcal H_\mu$. For details about the Verlinde formula we refer to 
\cite{huan21} and, for the case of WZW models, to \cite{beau8,sorg}. In the
special case $g \,{=}\, 0$ and $m \,{=}\, 3$, one obtains the statement that the 
matrix $S$ diagonalises the \emph{fusion rules} of the CFT, i.e.\ the multiplicity
matrices for the decompositions of the tensor products $\lambda_i \,{\otimes}\,
\lambda_j$ into simple objects. (The matrix $S$ also coincides with the matrix of
Hopf link invariants obtained in the three-dimensional surgery topological field
theory based on the modular tensor category $\mathrm{Rep}(\mathfrak V)$; that the 
latter matrix diagonalises the fusion rules follows quite directly 
\cite[Ch.\,IV.12]{TUra}, irrespectively of the connection with modular transformations.)


\subsubsection*{The relation with modular functors}

For a rational CFT there is a related more algebraic structure, also referred to 
as conformal blocks, in which the spaces of conformal blocks are
regarded just as finite-dimensional vector spaces 
that are endowed with a representation of the mapping class group, obey factorisation,
and depend functorially on the representation-theoretic data at the insertion points. 
These data are captured by the notion of a modular functor.
Various versions of this notion have been studied in the literature; 
roughly speaking, a modular functor is a three-dimensional topological field theory
that is only defined on a particular subclass of three-manifolds. For our purposes, and
in particular for the construction of CFT correlators further below, we work with the
following variant: an (anomaly-free) \emph{open-closed modular functor} is a symmetric 
monoidal pseudofunctor  
  \begin{equation}
  \mathrm{Bl}\colon\quad \mathrm{Bord}_{2,\mathrm{o/c}}^{\mathrm{or}}
  \xrightarrow{~~} \mathrm{Prof}_{\mathbb{C}}
  \label{eq:defmodularfunctor}
  \end{equation}
from the symmetric monoidal bicategory $\mathrm{Bord}_{2,\mathrm{o/c}}^{\mathrm{or}}$ 
to the symmetric monoidal bicategory $\mathrm{Prof}_{\mathbb{C}}$. The bicategory
$\mathrm{Bord}_{2,\mathrm{o/c}}^{\mathrm{or}}$ has compact oriented 1-manifolds as
objects, open-closed oriented two-bordisms as 1-morphisms, and isotopy classes of 
homeomorphisms (relative to boundary parametrisation) as 2-morphisms, while 
$\mathrm{Prof}_{\mathbb{C}}$ has (small) $\mathbb{C}$-linear categories as objects, 
$\mathbb{C}$-linear profunctors as 1-morphisms, and $\mathbb{C}$-linear natural 
transformations as 2-morphisms (a linear profunctor between linear categories $\mathcal A$ 
and $\mathcal B$ is a linear functor $P\colon \mathcal A^{\mathrm{op}}\,{\times}\, \mathcal B
\,{\xrightarrow{~}}\, \mathrm{Vect}_{\mathbb{C}}$.) The horizontal compositions
in $\mathrm{Prof}_{\mathbb{C}}$ are given by coends.

It is widely believed that any rational chiral CFT, with associated modular fusion
category $\mathcal C$, gives rise to a modular functor $\mathrm{Bl} \,{=}\, 
\mathrm{Bl}_{\mathcal C}$ such that the mapping class group representations are related 
by a Riemann-Hilbert correspondence to the monodromies of conformal blocks that come from
the projectively flat connection on the vector bundles $\mathcal B_{g,\vec\lambda}$.
Accordingly, yet another use of the term conformal blocks refers to the vector spaces
that are furnished by a modular functor.


\section*{Full conformal field theory}

As already stated, a full local conformal field theory is a consistent system of 
correlators  -- to be precise, of correlators for the class of world sheets one
is considering. The dependence on the conformal structure of a world sheet is
already incorporated in the modular functor \eqref{eq:defmodularfunctor}. Accordingly, 
from now on by a world sheet we mean a \emph{topological world sheet} $\mathcal{S}$, 
which has an underlying compact oriented topological 2-manifold $\varSigma_{\mathcal{S}}$
with possibly non-empty boundary. While in a chiral CFT the symmetries are encoded in
a vertex operator algebra $\mathfrak{V}$, in a full CFT we deal with the combined
action of two VOAs $\mathfrak{V}_{L}$ and $\mathfrak{V}_{R}$, which encode 
holomorphic and anti-holomorphic chiral symmetries, respectively, Accordingly the
relevant representation category is 
$\mathrm{Rep}(\mathfrak{V}_{L} \,{\otimes_{\mathbb C}}\,\mathfrak{V}_{R})$,
which under suitable finiteness conditions is equivalent to the Deligne product 
$\mathcal{C}_{L} \,{\boxtimes}\, \mathcal{C}_{R}$. This is referred to as the
\emph{combination of left and right movers}, or as \emph{holomorphic factorization}. In
the following we restrict our attention to rational CFTs for which the left and right 
movers are governed by the same VOA, such that $\mathcal{C}_{L} \,{=}\, \mathcal{C}$ and
$\mathcal{C}_{R} \,{=}\, \mathcal{C}^{\mathrm{rev}}$, where $\mathcal{C}^{\mathrm{rev}}$
equals $\mathcal{C}$ as a monoidal category, but has reversed braiding.
For a modular tensor category there is a canonical equivalence
  \begin{equation}
  \Xi\colon\quad \mathcal{C} \boxtimes \mathcal{C}^{\mathrm{rev}}
  \xrightarrow{~\simeq~} \mathcal{Z}(\mathcal{C})
  \label{eq:def:Xi}
  \end{equation}
of ribbon categories between the Deligne product 
$\mathcal{C} \,{\boxtimes}\, \mathcal{C}^{\mathrm{rev}}$ and the Drinfeld center
$\mathcal{Z}(\mathcal{C})$ of $\mathcal{C}$. (This relationship can be extended to
heterotic theories, for which $\mathcal{C}_{L}$ and $\mathcal{C}_{R}^{\mathrm{rev}}$ are
different \cite[Cor.\,3.30]{dmno}.)

We can then define a full CFT as a \emph{consistent system of correlators}, i.e.\
as an assignment 
  \begin{equation}
  \mathcal{S} \longmapsto \mathrm{Cor}_{\mathcal{C}}(\mathcal{S})
  \in\mathrm{Bl}_{\mathcal{C}}(\mathcal{S})
  \label{eq:def:conssystem}
  \end{equation}
that specifies for every world sheet $\mathcal{S}$ a correlator 
$\mathrm{Cor}_{\mathcal{C}}(\mathcal{S})$ as an element in the pertinent vector space
$\mathrm{Bl}_{\mathcal{C}}(\mathcal{S})$ of conformal blocks. Consistency of the system
means that $\mathrm{Cor}_{\mathcal{C}}(\mathcal{S})$ is required to be invariant under 
the action of the mapping class group $\mathrm{Map}(\mathcal{S})$ of the world sheet, 
and that the assignment is compatible with the sewing of world sheets. (In the presence 
of defects, $\mathrm{Map}(\mathcal{S})$ is typically a subgroup of the mapping class 
group $\mathrm{Map}(\varSigma_{\mathcal{S}})$ of the oriented topological surface 
$\varSigma_{\mathcal{S}}$.)
Given the category $\mathcal{C}$ of chiral data, the
problem of constructing a full CFT is of algebraic nature.

For oriented world sheets that are allowed to have physical boundaries, but no 
defects, a consistent system of correlators (subject to the extra requirements of
non-degeneracy of the sphere and disk two-point correlators and the uniqueness of
the closed state vacuum) is uniquely determined \cite{fjfrs,fjfrs2} by a 
\emph{simple special symmetric Frobenius algebra} $A$ internal to the modular fusion 
category $\mathcal{C}$, i.e.\ by a unital algebra object in $\mathcal{C}$ that has
non-zero quantum dimension, is simple as a bimodule over itself and is endowed with 
a counital coalgebra structure such that the comultiplication is a morphism of 
bimodules and is right inverse to the multiplication. (For unoriented world sheets,
the algebra $A$ must in addition be endowed with a so-called Jandl structure.)
The category $\mathrm{mod}\text{-}A$ of $A$-modules 
provides the conformal boundary conditions for such a CFT, and the algebra $A$ 
describes, up to Morita equivalence, the operator product expansions for boundary
fields that do not change the boundary condition. Equivalently, such an open-closed
full theory corresponds to an indecomposable semisimple $\mathcal{C}$-module
category $\mathcal{M}$ of boundary conditions. The mapping class
group invariance of a correlator reflects the fact that the corresponding
correlation function is a global section of the pertinent bundle of conformal blocks.

A larger class of world sheets is obtained when one also allows for topological defects. 
Defect conditions for topological defects of arbitrary codimension comprise a 
pivotal bicategory $\mathcal{F}r(\mathcal{C})$. The objects of $\mathcal{F}r(\mathcal{C})$
are simple symmetric special Frobenius algebras; they serve as labels for the 
\emph{phases} of the full CFT
which live on the two-dimensional strata of the world sheets. The 1- and 2-morphisms of 
$\mathcal{F}r(\mathcal{C})$ are bimodules and bimodules morphisms, which provide the
defect conditions for line and point defects, respectively. The pivotal bicategory
of defects can be equivalently realised as the bicategory 
$\mathcal{C}\text{-}\mathcal{M}\mathrm{od}^{\mathrm{tr}}$ of indecomposable semisimple 
$\mathcal{C}$-module categories admitting a module trace, module functors and 
module natural transformations.


\subsection*{Correlators from RT TFT and three-manifolds}

One way to construct a consistent system of correlators \cite{fuRs4,fjfrs} for a 
rational CFT is to exploit the connection with the three-dimensional topological 
field theory (TFT) of Reshetikhin-Turaev (RT) type that is associated to the modular 
fusion category $\mathcal{C}$. In this connection, the state spaces of the TFT 
functor $\mathrm{RT}_{\mathcal{C}}$ provide the spaces
of conformal blocks of the CFT: this is an instance of the holographic principle in 
quantum field theories and subsumes the duality between Chern-Simons TFTs and 
Wess-Zumino-Witten models as a special case. 

In this approach, called the \emph{TFT-construction}, one considers world sheets 
$\mathfrak{S}$ with field insertions. An illustrative example of such a world sheet is
  \begin{equation}
  \mathfrak{S}_0 ~:=~~ \scalebox{0.9}{\input{TFT0.tikz}}
  \label{eq:WSF}
  \end{equation}
This world sheet $\mathfrak{S}_0$ has the following attributes: two phases, labeled by 
Frobenius algebras $A$ and $B$, respectively, which we also indicate by using two 
different colors; three physical boundaries, with boundary conditions given by 
right $A$-modules $M$ and $M'$ and a right $B$-module $N$;
two line defects, with defect conditions given by $A$-$B$-bimodules $X$ and $Y$, 
and two point defects, with defect conditions $\alpha \,{\in}\,
\mathrm{Hom}_{\mathrm{mod}\text{-}B}(N,M\,{\otimes_{\!A}}\,X)$
and $\beta \,{\in}\, \mathrm{Hom}_{\mathrm{mod}\text{-}B}(M'\,{\otimes_{\!A}}\,Y,N)$.
Finally, $\mathfrak{S}_0$ has three field insertions: 

\noindent
(1) A \emph{boundary field} $\psi_{i}^{M,M'}$ which separates the two physical boundaries
labeled by $M$ and $M'$. It has one chiral label $i \,{\in}\, \mathrm{I}(\mathcal{C})$
given by a simple object in $\mathcal{C}$ together with an element of the degeneracy space 
  \begin{equation}
  \mathrm{Hom}_{\mathrm{mod}\text{-}A}(i\,{\otimes}\,M,M') \,.
  \label{eq:degspace4bdy}
  \end{equation}

\noindent
(2) A \emph{defect field} $\phi_{j,k}^{X,Y}$ which separates two line defects labeled 
by $X$ and $Y$. It has chiral labels $j,k \,{\in}\, \mathrm{I}(\mathcal{C})$ given by 
two simple objects in $\mathcal{C}$, combined with an element of the degeneracy space 
  \begin{equation}
  \mathrm{Hom}_{A\text{-}\mathrm{mod}\text{-}B}(j\,{\otimes^{-}}X\,{\otimes^{+}}k,Y) \,.
  \label{eq:degspace4bulk}
  \end{equation}
(The $A$-$B$-bimodule structure on $j\,{\otimes^{-}}X\,{\otimes^{+}}k$ is defined with
the help of the braiding of $\mathcal{C}$; for details see e.g.\ \cite{fuRs4}.)

\noindent
(3) A \emph{bulk field} $\varphi_{r,s}^{A}$ inserted in the region with phase $A$. 
It carries two chiral labels $r,s \,{\in}\, \mathrm{I}(\mathcal{C})$ and
an element of the degeneracy space 
$\mathrm{Hom}_{A\text{-}\mathrm{mod}\text{-}A}(r\,{\otimes^{-}}A\,{\otimes^{+}}s,A)$.

As suggested by the form of their degeneracy spaces, a bulk field is a special case of
a defect field: it separates two invisible defects, i.e.\ defect lines labeled by the 
Frobenius algebra (as a bimodule over itself). As an additional datum, the insertion
point $p$ of any field carries an arc-germ, i.e. an equivalence class of curves passing
through $p$ (two such curves are equivalent if they coincide in an open neighbourhood 
containing $p$). The arc-germ is what remains on the topological world sheet from
the germ of local coordinates on the original conformal world sheet.

 \medskip

In the TFT-construction, holomorphic factorization is implemented topologically
via the notion of the \emph{double} $\widehat{\varSigma}_{\mathfrak{S}}$ of a world sheet 
$\mathfrak{S}$: $\widehat{\varSigma}_{\mathfrak{S}}$ is a compact oriented topological 
surface with marked points (together with arc-germs) that is obtained as a
quotient of the orientation bundle over the surface $\varSigma_{\mathfrak{S}}$:
  \begin{equation}
  \widehat{\varSigma}_{\mathfrak{S}}\coloneqq\mathrm{Or}(\varSigma)/{\sim}
  \qquad\text{with}~~ (p,\mathrm{or})\,{\sim}\,(p,-\mathrm{or}) \quad \text{for}~
  p \,{\in}\, \partial\varSigma_{\mathfrak{S}} .
  \end{equation}
For instance, the double of the world sheet (\ref{eq:WSF}) is given
by a sphere with five marked points: %
  \begin{equation}
  \widehat{\varSigma}_{\mathfrak{S}_0} ~=\quad \input{TFT4.tikz}
  \end{equation}
The labels of the marked points come from the chiral labels of the field insertions.

The space of conformal blocks for a world sheet $\mathfrak{S}$ with field insertions is given
by evaluating the TFT functor on the double of $\mathfrak{S}$, i.e.\ $\mathrm{Bl}_{\mathcal{C}}
(\mathfrak{S}) \,{=}\, \mathrm{RT}_{\mathcal{C}}(\widehat{\varSigma}_{\mathfrak{S}})$.
In the TFT-construction, the correlator $\mathrm{Cor}_{\mathcal{C}}(\mathfrak{S})$ is 
obtained by evaluating $\mathrm{RT}_{\mathcal{C}}$ on a specific three-bordism 
$M_{\mathfrak{S}}\colon\emptyset \,{\xrightarrow{~}}\, \widehat{\varSigma}_{\mathfrak{S}}$
with an embedded $\mathcal{C}$-colored ribbon link; $\mathrm{Cor}_{\mathcal{C}}(\mathfrak{S})$
is the value at $1 \,{\in}\, \mathbb{C} \,{=}\, \mathrm{RT}_{\mathcal{C}}(\emptyset)$ of
the liner map $\mathrm{RT}_{\mathcal{C}}(M_{\mathfrak{S}})\colon \mathbb{C}
\,{\xrightarrow{~}}\, \mathrm{RT}_{\mathcal{C}}(\widehat{\varSigma}_{\mathfrak{S}})$.
As an illustration, the correlator of the world sheet \eqref{eq:WSF}
is obtained from the three-bordism 
  \begin{equation}
  M_{\mathfrak{S}_0} ~=\quad \scalebox{1.2}{\input{TFT1a.tikz}}
  \end{equation} 
The world sheet $\mathfrak{S}$ is a deformation retract of $M_{\mathfrak{S}}$; it can be 
seen sitting inside the solid ball. Each of its two-dimensional regions is replaced by a
network of ribbons located along a sufficiently fine dual triangulation of the region,
having trivalent vertices labeled by structural morphisms of the corresponding Frobenius
algebra. Each field insertion is replaced by its associated (bi)module morphism, with
the protruding legs colored with the respective chiral labels and attached to the 
marked points on $\partial M_{\mathfrak{S}} \,{=}\, \widehat{\varSigma}_{\mathfrak{S}}$.
For example, zooming in on the defect field insertion $\phi_{j,k}\equiv\phi_{j,k}^{X,Y}$
reveals the ribbon network
  \begin{equation}
  \input{TFT2.tikz}
  \end{equation}
(This explains the particular braidings appearing in \eqref{eq:degspace4bulk}.)

Making use of the defining properties of the bicategory $\mathcal{F}r(\mathcal{C})$, one
shows that the
correlator $\mathrm{Cor}_{\mathcal{C}}(\mathfrak{S}) \,{=}\, \mathrm{RT}_{\mathcal{C}}
(M_{\mathfrak{S}})(1) \,{\in}\, \mathrm{RT}_{\mathcal{C}}(\widehat{\varSigma}_{\mathfrak{S}})$
does not depend on the choice of triangulations and is
invariant under the mapping class group $\mathrm{Map}(\mathfrak{S})$
of the world sheet, which is the subgroup of $\mathrm{Map}(\varSigma_{\mathfrak{S}})$
containing those elements which fix the physical boundaries, defects
and field insertions up to isotopies. 

One insight resulting from the TFT-construction is that in the absence of field insertions
the correlators for a rational CFT can be expressed in terms of surface
defects that separate Reshetikhin-Turaev type TFTs \cite{kaSau3,fusV}.
When field insertions are present, there are in addition the ribbons carrying the chiral 
labels attached to the surface defect at the insertion points.
Also, instead of working with the doubles of world sheets, the holomorphic
factorization can alternatively be implemented by using the Reshetikhin-Turaev
TFT for the Drinfeld center $\mathcal{Z}(\mathcal{C})$.


\subsection*{Sewing boundaries and field contents}

When studying modular functors in
the framework of functorial field theory, it is convenient to describe a field not
via an insertion point carrying an arc-germ, but instead via a \emph{sewing boundary}.
Instead of world sheets $\mathfrak{S}$ with field insertions, we then deal with
world sheets $\mathcal{S}$ with sewing boundaries. For instance, the world sheet 
$\mathfrak{S}_0$ in \eqref{eq:WSF} gets replaced by the world sheet 
  \begin{equation}
  \mathcal{S}_0 ~=\quad \scalebox{0.9}{\input{TFT3.tikz}}
  \label{eq:WSS}
  \end{equation}
$\mathcal{S}_0$ has one \emph{sewing interval} which separates two physical boundaries
and two \emph{sewing circles} obtained from cutting out disks around the defect and 
bulk field insertions.  World sheets can be sewn along sewing boundaries with matching 
boundary data. For every sewing boundary there is a field content, i.e.\ a space 
of fields associated to it. For instance, the field for the sewing boundaries of the
world sheet \eqref{eq:WSS} are:

\noindent
(1) The boundary field content 
$\mathbb{B}^{M,M'}\,{:=}\,\bigoplus_{i\in\mathrm{I}(\mathcal{C})}
\mathrm{Hom}_{\mathrm{mod}\text{-}A}(i\,{\otimes}\, M,M')\otimes_{\mathbb{C}}i$\,.

\noindent
(2) The defect field content 
$\mathbb{D}^{X,Y}\,{:=}\,\bigoplus_{i,j\in\mathrm{I}(\mathcal{C})}
\mathrm{Hom}_{A\text{-}\mathrm{mod}\text{-}B}(i\,{\otimes^{-}}X\,{\otimes^{+}}j,Y)
\otimes_{\mathbb{C}} \Xi(i\,{\boxtimes}\, j)$\,. 

\noindent
(3) The bulk field content 
$\mathbb{D}^{A,A}\,{:=}\,\bigoplus_{i,j\in\mathrm{I}(\mathcal{C})}
\mathrm{Hom}_{A\text{-}\mathrm{mod}\text{-}A}(i\,{\otimes^{-}}A\,{\otimes^{+}}j,A)
\otimes_{\mathbb{C}} \Xi(i\,{\boxtimes}\, j)$\,. 

Here $\Xi$ is the ribbon equivalence \eqref{eq:def:Xi} between 
$\mathcal{C} \,{\boxtimes}\, \mathcal{C}^{\mathrm{rev}}$ and the Drinfeld center
$\mathcal{Z}(\mathcal{C})$.
A boundary field content is an object in $\mathcal{C} \,{=}\, \mathrm{Rep}(\mathfrak{V})$
because we require the boundary conditions to preserve the chiral VOA $\mathfrak{V}$.
In contrast, while the spaces of defect and bulk fields are
naturally modules over $\mathfrak{V}\,{\otimes_{\mathbb{C}}}\,\mathfrak{V}$, hence
their field contents are objects in
$\mathcal{Z}(\mathcal{C})\,{\simeq}\,\mathcal{C} \,{\boxtimes}\, \mathcal{C}^{\mathrm{rev}}$.

The field contents can also be expressed as internal Homs. For a left module
category $\mathcal{M}$ over a monoidal category $\mathcal{A}$ and a pair of objects 
$m,n \,{\in}\, \mathcal{M}$, the internal Hom $\underline{\mathrm{Hom}}_{\mathcal{M}}(m,n)$
(if it exists, which is the case especially when $\mathcal{M}$ and $\mathcal{A}$
are finitely semisimple) is an object in $\mathcal{A}$ that is defined up to unique 
isomorphisms by the adjunction $\mathrm{Hom}_{\mathcal{M}}(c\,{\triangleright}\, m,n)
\,{\cong}\, \mathrm{Hom}_{\mathcal{A}}(c,\underline{\mathrm{Hom}}_{\mathcal{M}}(m,n))$.
Since $\mathrm{mod}\text{-}A$ is a finitely semisimple left module
category over the fusion category $\mathcal{C}$, we have
  \begin{equation}
  \begin{aligned}
  \mathbb{B}^{M,M'} & =\bigoplus_{i\in\mathrm{I}(\mathcal{C})}
  \mathrm{Hom}_{\mathrm{mod}\text{-}A}(i\,{\otimes}\, M,M')\otimes_{\mathbb{C}}i
  = \int^{c\in\mathcal{C}}\!\! \mathrm{Hom}_{\mathrm{mod}\text{-}A}
  (c\,{\triangleright}\, M,M') \otimes_{\mathbb{C}} c
  \\ &
  \cong \int^{c\in\mathcal{C}}\! \mathrm{Hom}_{\mathcal{C}}
  (c,\underline{\mathrm{Hom}}_{\mathrm{mod}\text{-}A}(M,M')) \otimes_{\mathbb{C}} c
  = \underline{\mathrm{Hom}}_{\mathrm{mod}\text{-}A}(M,M') \,\in \mathcal{C} \,.
  \end{aligned}
  \end{equation}
Here the integral sign denotes a categorical coend. If $M \,{=}\, M'$ and the quantum 
dimension $\mathrm{dim}(M)$ (of $M$ as an object in $\mathcal{C}$) is non-zero, then
$\mathbb{B}^{M,M} \,{=}\, \underline{\mathrm{Hom}}_{\mathrm{mod}\text{-}A}(M,M)$
is a simple special symmetric Frobenius algebra
that is Morita equivalent to $A$. By a similar, albeit less straightforward, procedure 
for the finite left $\mathcal{Z}(\mathcal{C})$-module category $\mathcal Fun_{A,B}$
of $\mathcal{C}$-module
functors from $\mathrm{mod}\text{-}A$ to $\mathrm{mod}\text{-}B$ it follows
\cite{fuSc25} that defect field contents are \emph{internal natural transformations}
  \begin{equation}
  \mathbb{D}^{X,Y} = \underline{\mathrm{Nat}}(-\,{\otimes_{A}}\,X,-\,{\otimes_{A}}\,Y)
  := \underline{\mathrm{Hom}}_{\mathcal Fun_{A,B}}(-\,{\otimes_{A}}\,X,-\,{\otimes_{A}}\,Y)
  \,\in \mathcal{Z}(\mathcal{C}) \,.
  \end{equation}


The bulk field content $\mathbb{D}^{A,A} \,{=}\, \underline{\mathrm{Nat}}
(\mathrm{id}_{\mathrm{mod}\text{-}A},\mathrm{id}_{\mathrm{mod}\text{-}A})
\,{=:}\, Z(A) \,{\in}\, \mathcal{Z}(\mathcal{C})$,
called the full center of the simple special symmetric Frobenius
algebra $A \,{\in}\, \mathcal{C}$, is a commutative symmetric
Frobenius algebra in the Drinfeld center $\mathcal{Z}(\mathcal{C})$;
$Z(A)$ is \emph{Lagrangian} in the sense of \cite[Def.\,4.6]{dmno}.


\subsection*{Modular functor, field maps and correlators}

We now formulate the full local rational CFT with given modular fusion category
$\mathcal{C}$ of chiral data with the help of an open-closed modular functor 
$\mathrm{Bl}_{\mathcal{C}}$, as defined in \eqref{eq:defmodularfunctor}. 
$\mathrm{Bl}_{\mathcal{C}}$ provides the conformal blocks of the CFT. For
obtaining correlators, $\mathrm{Bl}_{\mathcal{C}}$ must
in addition fulfill the following requirements: 
   \begin{itemize}
\item 
The categories assigned to the interval $I \,{=}\, [0,1]\,{\subset}\,\mathbb{R}$ and to
the circle $S^{1} \,{=}\, \{|z| \,{=}\, 1\} \,{\subset}\, \mathbb{C}$ are equipped with 
equivalences
  \begin{equation}
  \Phi_{I}\colon~~ \mathrm{Bl}_{\mathcal{C}}(I)\xrightarrow{~\simeq~} \mathcal{C}
  \qquad \text{and} \qquad 
  \Phi_{S^{1}}\colon~~ \mathrm{Bl}_{\mathcal{C}}(S^{1})\xrightarrow{~\simeq~}
  \mathcal{Z}(\mathcal{C})
  \label{eq:PhiI.PhiS}
  \end{equation}
of $\mathbb{C}$-linear categories. 
\item 
The closed sector of the modular functor $\mathrm{Bl}_{\mathcal{C}}$
is canonically equivalent to the Reshetikhin-Turaev modular functor
$\mathrm{RT}_{\!\mathcal{Z}(\mathcal{C})}$.
\end{itemize}
The equivalences \eqref{eq:PhiI.PhiS} extend uniquely to equivalences 
$\Phi_{\ell}\colon\mathrm{Bl}_{\mathcal{C}}(\ell)\,{\xrightarrow{\,\simeq\,}}\, 
\mathcal{C}^{\times p} \,{\times}\, \mathcal{Z}(\mathcal{C})^{\times q}$
for all one-manifolds $\ell \,{=}\, (I)^{\sqcup p} \,{\sqcup}\, (S^{1})^{\sqcup q}$.

A second ingredient needed for the construction of correlators is a collection of 
\emph{field maps} which encode the field contents of all types of fields, including 
multi-pronged generalizations of defect fields, which are not considered traditionally. To
describe these we need the notion of an \emph{$\mathcal{F}r(\mathcal{C})$-boundary datum}
$\mathsf{b}$ on a compact oriented one-manifold $\ell$; this is a finite set 
$O \,{\subset}\, \mathrm{int}(\ell)$ of points in the interior of $\ell$, together with 
a labeling of the connected components of the complement $\ell\,{\setminus}\, O$ by
the objects in $\mathcal{F}r(\mathcal{C})$, i.e.\ by simple special symmetric Frobenius
algebras, and a labeling of the elements of $O$ by the 1-morphisms in 
$\mathcal{F}r(\mathcal{C})$, i.e.\ by bimodules.
Denote by $\mathcal{F}r(\mathcal{C})_{\ell}$ the set of $\mathcal{F}r(\mathcal{C})$-boundary
data on $\ell$. Given a world sheet $\mathcal{S}$, a structure of an open-closed two-bordism
$\varSigma_{\mathcal{S}}\colon\ell_{\mathrm{in}} \,{\xrightarrow{~}}\, \ell_{\mathrm{out}}$
on its underlying surface uniquely determines $\mathcal{F}r(\mathcal{C})$-boundary
data $\mathsf{b}_{\mathrm{in}} \,{\in}\, \mathcal{F}r(\mathcal{C})_{\ell_{\mathrm{in}}}$
and $\mathsf{b}_{\mathrm{out}} \,{\in}\, \mathcal{F}r(\mathcal{C})_{\ell_{\mathrm{out}}}$.
Field maps are a collection $\{\mathbb{F}_{\ell}\colon\mathcal{F}r(\mathcal{C})_{\ell}
\,{\xrightarrow{~}}\, \mathrm{obj}(\mathrm{Bl}_{\mathcal{C}}(\ell))\}_{\ell}$
of maps defined for every compact oriented one-manifold $\ell$ with any numbers $p$ 
of intervals and $q$ of circles as connected components, such that for any world
sheet $\mathcal{S}$ with its underlying surface viewed as an open-closed bordism, 
the objects $\Phi_{\ell_{\varepsilon}} \,{\circ}\, \mathbb{F}_{\ell_{\varepsilon}}
(\mathsf{b}_{\varepsilon}) \,{\in}\,\mathcal{C}^{\times p_{\varepsilon}}
{\times}\, \mathcal{Z}(\mathcal{C})^{\times q_{\varepsilon}}$, for
$\varepsilon \,{\in}\, \{\mathrm{in},\mathrm{out}\}$, are given by the correct
field contents.

Given the open-closed modular functor $\mathrm{Bl_{\mathcal{C}}}$ and the field maps
$\{\mathbb{F}_{\ell}\}$, we obtain the vector spaces of conformal blocks as follows. The
space $\mathrm{Bl}_{\mathcal{C}}(\mathcal{S})$
of conformal blocks for a world sheet $\mathcal{S}$ with underlying bordism
$\varSigma_{\mathcal{S}}\colon\ell_{\mathrm{in}} \,{\xrightarrow{~}}\, \ell_{\mathrm{out}}$
is the vector space
  \begin{equation}
  \mathrm{Bl}_{\mathcal{C}}(\mathcal{S})\,{:=}\, \mathrm{Bl}_{\mathcal{C}}
  (\varSigma_{\mathcal{S}};\mathbb{F}_{\ell_{\mathrm{in}}}(\mathsf{b}_{\mathrm{in}}),
  \mathbb{F}_{\ell_{\mathrm{out}}}(\mathsf{b}_{\mathrm{out}})) \,,
  \end{equation}
where $\mathrm{Bl}_{\mathcal{C}}(\varSigma_{\mathcal{S}},-;\sim)\colon
\mathrm{Bl}_{\mathcal{C}}(\ell_{\mathrm{in}})^{\mathrm{op}} \,{\times}\,
\mathrm{Bl}_{\mathcal{C}}(\ell_{\mathrm{out}}) \,{\xrightarrow{~}}\, \mathrm{Vect}_{\mathbb{C}}$
is the profunctor obtained by evaluating the modular functor on $\varSigma_{\mathcal{S}}$.
By functoriality, the vector space $\mathrm{Bl}_{\mathcal{C}}(\mathcal{S})$
carries an action of the mapping class group $\mathrm{Map}(\varSigma_{\mathcal{S}}) \,{=}\,
\mathrm{End}_{\mathrm{Bord}_{2,\mathrm{o/c}}^{\mathrm{or}}}(\varSigma_{\mathcal{S}})$.
Also by functoriality, for any sewing of two world sheets $\mathcal{S}$ and $\mathcal{S}'$
along a sewing boundary $\ell$ we get a \emph{sewing map} 
  \begin{equation}
  s\colon\quad \mathrm{Bl}_{\mathcal{C}}(\mathcal{S}) \otimes_{\mathbb{C}}
  \mathrm{Bl}_{\mathcal{C}}(\mathcal{S}')
  \,{\xrightarrow{~}}\, \mathrm{Bl}_{\mathcal{C}}(\mathcal{S}\cup_{\ell}\mathcal{S}') 
  \label{eq:def:sewingmap}
  \end{equation}
of conformal blocks.
These maps endow $\mathrm{Bl}_{\mathcal{C}}(\mathcal{S}\cup_{\ell}\mathcal{S}')$ 
with the structure of a coend, 
  \begin{equation}
  \mathrm{Bl}_{\mathcal{C}}(\mathcal{S}{\cup_{\ell}}\mathcal{S}')
  = \int^{\mathsf{b}\in\mathrm{Bl}_{\mathcal{C}}(\ell)}\! \mathrm{Bl}_{\mathcal{C}}
  (\varSigma_{\mathcal{S}};\mathbb{F}_{\ell_{\mathrm{in}}}(\mathsf{b}_{\mathrm{in}}),
  \mathsf{b}) \,{\otimes_{\mathbb{C}}}\, \mathrm{Bl}_{\mathcal{C}}(\varSigma_{\mathcal{S}};
  \mathsf{b}_{\mathrm{}},\mathbb{F}_{\ell'_{\mathrm{out}}}(\mathsf{b}'_{\mathrm{out}})) \,.
  \label{eq:coendBlC}
  \end{equation}
A consistent system of correlators is then an assignment 
$\mathcal{S} \,{\mapsto}\, \mathrm{Cor}_{\mathcal{C}}(\mathcal{S})
\,{\in}\,\mathrm{Bl}_{\mathcal{C}}(\mathcal{S})$, as in \eqref{eq:def:conssystem},
such that 
$\mathrm{Cor}_{\mathcal{C}}(\mathcal{S})$ is invariant under the action of 
$\mathrm{Map}(\mathcal{S}) \,{\subset}\, \mathrm{Map}(\varSigma_{\mathcal{S}})$
and such that the sewing maps \eqref{eq:def:sewingmap} take correlators to correlators.


\subsection*{Correlators from string nets}

Another approach to the construction of correlators, called the 
\emph{string-net construction} \cite{fusY}, uses a realization of the open-closed 
modular functor $\mathrm{Bl}_{\mathcal{C}}$ via string nets. String-net
models arose in the study of topologically ordered phases of matter \cite{leWe} and 
were later formulated as two-dimensional skein theories \cite{kirI24}.

The basic input datum of a string-net model is a spherical fusion category. In our
context this is the modular fusion category $\mathcal{C}$ of chiral data. For any 
compact oriented surface $\varSigma$ and $\mathcal{C}$-boundary datum 
$\mathsf{B}^{\circ}$ on $\partial\varSigma$ (defined in the same way as an 
$\mathcal{F}r(\mathcal{C})$-boundary datum, by regarding $\mathcal{C}$ as a one-object
bicategory), the string-net model for $\mathcal{C}$ defines a finite-dimensional
vector space $\mathrm{SN}_{\mathcal{C}}^{\circ}(\varSigma,\mathsf{B}^{\circ})$,
called the bare string-net space associated to $(\varSigma,\mathsf{B}^{\circ})$,
as a quotient of the free vector space generated by $\mathcal{C}$-colored
string diagrams drawn on the surface $\varSigma$ with boundary datum
$\mathsf{B}^{\circ}$, by a subspace that encodes the local graphical calculus of 
$\mathcal{C}$. Thus a \emph{string net}, i.e.\ an element of 
$\mathrm{SN}_{\mathcal{C}}^{\circ}(\varSigma,\mathsf{B}^{\circ})$,
is a linear combination of equivalence classes of string diagrams on $\varSigma$, 
where two string diagrams are equivalent if they can be transformed into each 
other by applying the graphical calculus of $\mathcal{C}$ within disk-shaped regions.

For every compact oriented one-manifold $\ell$ one can then define a linear
category $\mathrm{Cyl^{\circ}}(\mathcal{C},\ell)$, called the bare cylinder category
over $\ell$, whose objects are $\mathcal{C}$-boundary data on $\ell$ and whose 
morphisms are string nets on the cylinder $\ell\,{\times}\, I$. Composition in 
$\mathrm{Cyl^{\circ}}(\mathcal{C},\ell)$ is given by sewing the cylinders and 
concatenating the string diagrams. The \emph{cylinder category}
$\mathrm{Cyl}(\mathcal{C},\ell)$ over $\ell$ is the idempotent completion of 
$\mathrm{Cyl^{\circ}}(\mathcal{C},\ell)$, whose objects are idempotents in 
$\mathrm{Cyl^{\circ}}(\mathcal{C},\ell)$. We have canonical equivalences 
$\mathrm{Cyl}(\mathcal{C},\ell) \,{\xrightarrow{\,\simeq\,}}\, \mathcal{C}$ and 
$\mathrm{Cyl}(\mathcal{C},S^{1}) \,{\xrightarrow{\,\simeq\,}}\, \mathcal{Z}(\mathcal{C})$
of linear categories. Accordingly one defines the string-net space 
$\mathrm{SN}_{\mathcal{C}}(\varSigma,\mathsf{B})$, which takes an object $\mathsf{B}$ 
in $\mathrm{Cyl}(\mathcal{C},\ell)$ as its boundary datum, as the subspace of the 
bare string-net space that consists of elements which are invariant
under sewing with the string net on a cylinder that is given by $\mathsf{B}$.
The so defined string-net spaces carry a mapping class group action
obtained by pushforward. When $\varSigma$ is equipped with the structure
of an open-closed bordism, one obtains a profunctor 
  \begin{equation}
  \mathrm{SN}_{\mathcal{C}}(\varSigma)\colon\quad
  \mathrm{Cyl}(\mathcal{C},\ell_{\mathrm{in}})^{\mathrm{op}} \times 
  \mathrm{Cyl}(\mathcal{C},\ell_{\mathrm{out}}) \xrightarrow{~~} \mathrm{Vect} \,.
  \end{equation}
The assignments $\ell\,{\mapsto}\,\mathrm{Cyl}(\mathcal{C},\ell)$ and
$\varSigma\,{\mapsto}\,\mathrm{SN}_{\mathcal{C}}(\varSigma)$ define an open-closed
modular functor $\mathrm{SN}_{\mathcal{C}}$. Moreover, the closed sector of 
$\mathrm{SN}_{\mathcal{C}}$ extends to a once-extended three-dimensional TFT 
that is equivalent to the Turaev-Viro TFT $\mathrm{TV}_{\mathcal{C}}$
which, in turn, is equivalent 
to $\mathrm{RT}_{\mathcal{Z}(\mathcal{C})}$.
Thus one can use the string-net modular functor $\mathrm{SN}_{\mathcal{C}}$
as the model for conformal blocks.

In the string-net approach, the construction of correlators is fairly straightforward. The
correlator $\mathrm{Cor}_{\mathcal{C}}(\mathcal{S})$ is a $\mathcal{C}$-colored string
net that is represented by a string diagram $\varGamma_{\!\mathcal{C}}(\mathcal{S})$ 
that is obtained by replacing each two-dimensional stratum of $\mathcal{S}$
with a network of strings labeled by the relevant Frobenius algebra according to a 
fine triangulation,
and, when $\mathcal{S}$ has physical boundaries, by adding 
two-dimensional strata labeled by the trivial Frobenius algebra $\mathbf1$ which turn
a physical boundary (a right $A$-module) into a defect line (a $\mathbf1$-$A$-bimodule).
As an illustration, a string diagram for the world sheet
  \begin{equation}
  \mathcal{S}_1 ~= \scalebox{1.4}{\input{WS1.tikz}}
  \end{equation} 
is given by
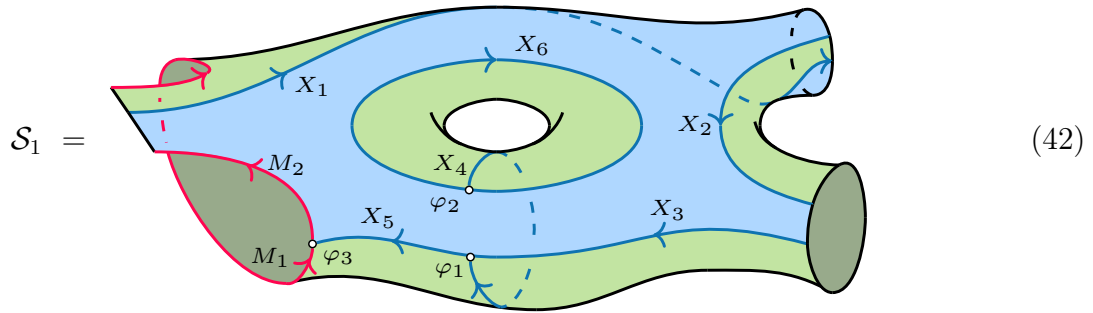
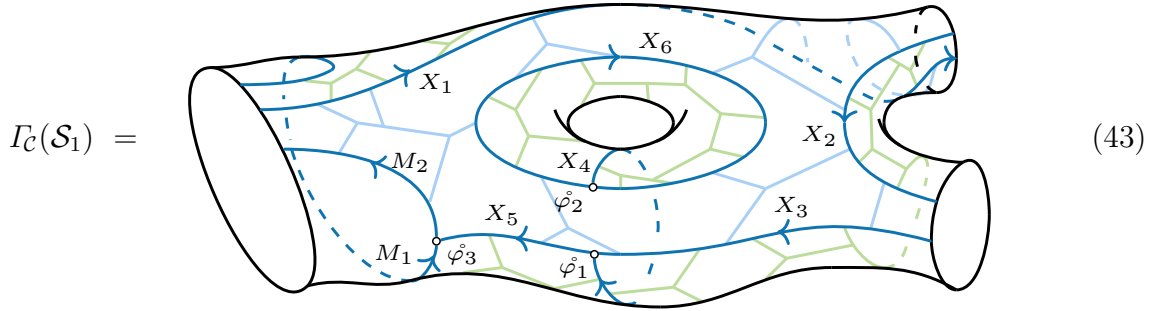
  \begin{equation}
  \varGamma_{\!\mathcal{C}}(\mathcal{S}_1) ~=\!\!\!\! \scalebox{1.4}{\input{WS3.tikz}}
  \end{equation}

Due to the defining properties of the underlying algebraic structures,
the so obtained string net is well defined and is invariant under the action of 
$\mathrm{Map}(\mathcal{S})$. Also, it is in the string-net space 
$\mathrm{SN}_{\mathcal{C}}(\varSigma_{\mathcal{S}},
\mathbb{F}_{\!\partial\varSigma_{\mathcal{S}}}(\mathsf{b}_{\mathcal{S}}))$,
which is taken to be the space $\mathrm{Bl_{\mathcal{C}}(\mathcal{S})}$ of conformal
blocks, where the boundary datum 
$\mathbb{F}_{\!\partial\varSigma_{\mathcal{S}}}(\mathsf{b}_{\mathcal{S}}) \,{\in}\, 
\mathrm{Cyl}(\mathcal{C},\partial\varSigma_{\mathcal{S}}) \,{\simeq}\,
\mathcal{C}^{\times p} \,{\times}\, \mathcal{Z}(\mathcal{C})^{\times q}$
produces the correct field contents.

Upon choosing a marking without cuts (in the sense of \cite{baKir}) for the surface, 
a correlator can be translated into a morphism in $\mathcal{Z}(\mathcal{C})$. The 
string-net description of correlators for world sheets of particular interest, such 
as those that correspond to operator products, provides explicit expressions which
match results proposed in the literature. For instance, operator products of defect 
fields correspond to the vertical and horizontal compositions of internal natural 
transformations, and in the special case of a bulk field they reduce to the 
commutative product of the full center. Moreover, the torus partition function 
$\mathrm{Cor}_{\mathcal C}(\mathcal T_{\!A})\,{\in}\,\mathrm{SN}_{\mathcal C}(\mathrm T)$,
where the world sheet $\mathcal{T}_{\!A}$ is a torus $\mathrm T$ without defect lines
whose phase is given by a Frobenius algebra $A \,{\in}\, \mathcal{F}r(\mathcal{C})$, 
decomposes as $\mathrm{Cor}_{\mathcal{C}}(\mathcal{T}_{A}) \,{=}\, 
\sum_{i,j\in\mathrm{I}(\mathcal{C})}Z_{i,j}(A)\,e_{i,j}$, with the integers
  \begin{equation}
  Z_{i,j}(A) = \mathrm{dim}_{\mathbb{C}}(
  \mathrm{Hom}_{A\text{-}\mathrm{mod}\text{-}A}(i\,{\otimes^{-}}A\,{\otimes^{+}}j,A))
  \end{equation}
being dimensions of degeneracy spaces for bulk fields (compare \eqref{eq:degspace4bulk}).
Here $\{e_{i,j}\}_{i,j\in\mathrm{I}(\mathcal{C})}$ is a distinguished basis of 
$\mathrm{SN}_{\mathcal{C}}(\mathrm{T})$ 
that corresponds to the characters \eqref{eq:def:char} of VOA-modules
and transforms under the modular group $\mathrm{Map}(\mathrm T) \,{=}\, SL(2,\mathbb{Z})$
via conjugation by the modular $S$- and $T$-matrices of $\mathcal{C}$. The mapping 
class group invariance of the string net $\mathrm{Cor}_{\mathcal{C}}(\mathcal{T}_{\!A})$
immediately implies the well known modular invariance relations
$[Z(A),S] \,{=}\, 0 \,{=}\, [Z(A),T]$.

 \medskip

For the implementation of categorical symmetries \cite{ffrs5,frmT} via topological 
defects, a desirable property of the system of correlators is that one can modify the
defect network in accordance with the local graphical calculus for the category
$\mathcal{C}$ on any world sheet 
without altering its correlator. The string-net correlators indeed enjoy this property;
this can be shown with the help of generalised string-net models which take,
instead of a spherical fusion category, a \emph{pivotal bicategory} as input datum
\cite{fusY2}. Taking this input bicategory to be $\mathcal{F}r(\mathcal{C})$, one
defines for any compact oriented surface $\varSigma$ and
$\mathcal{F}r(\mathcal{C})$-boundary datum $\mathsf{b}$ a vector space 
$\mathrm{SN}_{\mathcal{F}r(\mathcal{C})}^{\circ}(\varSigma,\mathsf{b})$
of $\mathcal{F}r(\mathcal{C})$-colored string nets on $\varSigma$. The string-net 
space $\mathrm{SN}_{\mathcal{F}r(\mathcal{C})}^{\circ}(\varSigma,\mathsf{b})$ can be
interpreted as the space of equivalence classes of world sheets with underlying surface
$\varSigma$ and boundary datum $\mathsf{b}$, with the equivalence relation given by 
the local graphical calculus of the pivotal bicategory $\mathcal{F}r(\mathcal{C})$. It
is thus appropriate to call an element of the space
$\mathrm{SN}_{\mathcal{F}r(\mathcal{C})}^{\circ}(\varSigma,\mathsf{b})$
a \emph{quantum world sheet}. The assignment of the correlators provides a linear map
$\mathrm{Cor}_{\mathcal{C}}(\varSigma,\mathsf{b}) \colon
\mathbb{C}\mathrm{G}_{\mathcal{F}r(\mathcal{C})}(\varSigma,\mathsf{b}) \,{\xrightarrow{~}
}\, \mathrm{SN}_{\mathcal{C}}(\varSigma,\mathbb{F}_{\!\partial\varSigma}(\mathsf{b}))$,
where $\mathbb{C}\mathrm{G}_{\mathcal{F}r(\mathcal{C})}(\varSigma,\mathsf{b})$
is the vector space generated by the set of $\mathcal{F}r(\mathcal{C})$-colored
string diagrams on $\varSigma$. This map factors through the canonical quotient map 
$\mathrm{q}(\varSigma,\mathsf{b}) \colon \mathbb{C}\mathrm{G}_{\mathcal{F}r(\mathcal{C})}
(\varSigma,\mathsf{b}) \,{\twoheadrightarrow}\, \mathrm{SN}_{\mathcal{F}r
(\mathcal{C})}^{\circ}(\varSigma,\mathsf{b})$, i.e.\ there is a unique linear map 
  \begin{equation}
  \mathrm{U}(\varSigma,\mathsf{b})\colon\quad
  \mathrm{SN}_{\mathcal{F}r(\mathcal{C})}^{\circ}(\varSigma,\mathsf{b}) \xrightarrow{~~}
  \mathrm{SN}_{\mathcal{C}}(\varSigma,\mathbb{F}_{\!\partial\varSigma}(\mathsf{b}))
  \end{equation}
such that $\mathrm{Cor}_{\mathcal{C}}(\varSigma,\mathsf{b}) \,{=}\,
\mathrm{U}(\varSigma,\mathsf{b}) \,{\circ}\, \mathrm{q}(\varSigma,\mathsf{b})$.
As a consequence, the correlator of a world sheet only depends on the
quantum world sheet it represents, and is therefore unchanged under modifications
via the local graphical calculus. The map $\mathrm{U}(\varSigma,\mathsf{b})$ is
called a \emph{universal correlator}. It intertwines with the action
of the mapping class group.


\section*{Conclusion}

We have presented various algebraic structures that play a role in 
conformal field theory, considering both chiral and full local CFTs. It
is worth pointing out that all of these are perfectly customary mathematical
structures and that they can be analysed with standard mathematical tools.

Needless to say, our exposition is highly biased by the authors' taste and restricted
knowledge. Many interesting aspects and important developments
as well as a substantial portion of pertinent literature had to be omitted owing 
to limitations of length and of expertise of the authors.


\vskip 3em

\noindent
{\sc Acknowledgements:}\\[.3em]
We thank Chiara Damiolini for helpful comments on an earlier version of the text.
JF is supported by VR under project no.\ 2022-02931. CS is supported by 
the Deutsche Forschungsgemeinschaft (DFG, German Research Foundation) under
SCHW1162/6-1 and under Germany's Excellence Strategy - EXC 2121 ``Quantum Universe'' 
- 390833306.
SW is supported by the Engineering and Physical Sciences Research Council (EPSRC)
EP/V053787/1 and by the Alexander von Humboldt Foundation.
YY is supported by a Junior Research Fellowship of the
Erwin Schr\"odinger International Institute for Mathematics and Physics (ESI).

\newpage
\pagestyle{plain}

\end{document}

%% file: sample.tikzstyles
\tikzstyle{nc}=[fill=white, draw=black, shape=circle, minimum size=0.40cm]
\tikzstyle{nr}=[fill=white, draw=black, shape=rectangle, minimum width=0.60cm, minimum height=0.60cm]
\tikzstyle{nrg}=[fill={rgb,255: red,199; green,199; blue,199}, draw=black, shape=rectangle]
\tikzstyle{sr}=[fill=white, draw=black, shape=rectangle, minimum width=0.60cm, minimum height=0.40cm, inner sep=0pt]
\tikzstyle{lsr}=[fill=white, draw=black, shape=rectangle, minimum width=0.60cm, minimum height=0.40cm, inner sep=0pt, rotate=30]
\tikzstyle{rsr}=[fill=white, draw=black, shape=rectangle, minimum width=0.60cm, minimum height=0.40cm, inner sep=0pt, rotate=-30]
\tikzstyle{rrrsr}=[fill=white, draw=black, shape=rectangle, minimum width=0.60cm, minimum height=0.40cm, inner sep=0pt, rotate=-120]
\tikzstyle{ivsr}=[fill=none, draw=none, shape=rectangle, minimum width=0.60cm, minimum height=0.40cm, inner sep=0pt]
\tikzstyle{ivlsr}=[fill=none, draw=none, shape=rectangle, minimum width=0.60cm, minimum height=0.40cm, inner sep=0pt, rotate=30]
\tikzstyle{ivrsr}=[fill=none, draw=none, shape=rectangle, minimum width=0.60cm, minimum height=0.40cm, inner sep=0pt, rotate=-30]
\tikzstyle{ivrrrsr}=[fill=none, draw=none, shape=rectangle, minimum width=0.60cm, minimum height=0.40cm, inner sep=0pt, rotate=-120]
\tikzstyle{mr}=[fill=white, draw=black, shape=rectangle, minimum width=1.00cm, minimum height=0.50cm, inner sep=0pt]
\tikzstyle{product}=[fill={rgb,255: red,255; green,5; blue,80}, draw=black, shape=circle, inner sep=0pt, minimum size=0.15cm]
\tikzstyle{unit}=[fill={rgb,255: red,255; green,166; blue,217}, draw=black, shape=circle, inner sep=0pt, minimum size=0.15cm]
\tikzstyle{coproduct}=[fill={rgb,255: red,2; green,145; blue,255}, draw=black, shape=circle, inner sep=0pt, minimum size=0.15cm]
\tikzstyle{counit}=[fill={rgb,255: red,165; green,219; blue,255}, draw=black, shape=circle, inner sep=0pt, minimum size=0.15cm]
\tikzstyle{antipode}=[fill=white, draw=black, shape=circle, inner sep=0pt, minimum size=0.15cm]
\tikzstyle{region}=[fill={rgb,255: red,191; green,191; blue,191}, draw={rgb,255: red,191; green,191; blue,191}, shape=circle, minimum size=0.80cm]
\tikzstyle{boundarydisc}=[fill={rgb,255: red,128; green,128; blue,128}, draw=black, shape=circle, minimum size=0.80cm]
\tikzstyle{identity}=[fill=black, draw=black, shape=circle, inner sep=0pt, minimum size=0.15cm]
\tikzstyle{S}=[fill=white, draw=black, shape=rectangle, minimum width=0.15cm, minimum height=0.15cm]
\tikzstyle{S-}=[fill={rgb,255: red,199; green,199; blue,199}, draw=black, shape=rectangle, minimum width=0.15cm, minimum height=0.15cm]
\tikzstyle{dot}=[fill=white, draw=black, shape=circle, inner sep=0pt, minimum size=0.07cm]
\tikzstyle{bdot}=[fill=black, draw=black, shape=circle, inner sep=0pt, minimum size=0.07cm]
\tikzstyle{m}=[fill={rgb,255: red,255; green,5; blue,80}, draw={rgb,255: red,255; green,5; blue,80}, shape=circle, inner sep=0pt, minimum size=0.14cm, tikzit shape=circle]
\tikzstyle{mod}=[fill={rgb,255: red,172; green,229; blue,232}, draw={rgb,255: red,14; green,115; blue,177}, shape=circle]
\tikzstyle{smod}=[fill={rgb,255: red,172; green,229; blue,232}, draw={rgb,255: red,14; green,115; blue,177}, shape=circle, inner sep=0pt, minimum size=0.15cm]
\tikzstyle{smod1}=[fill={rgb,255: red,150; green,255; blue,138}, draw={rgb,255: red,5; green,130; blue,3}, shape=circle, inner sep=0pt, minimum size=0.15cm]
\tikzstyle{smod2}=[fill={rgb,255: red,255; green,202; blue,239}, draw={rgb,255: red,126; green,0; blue,124}, shape=circle, inner sep=0pt, minimum size=0.15cm]
\tikzstyle{smod3}=[fill={rgb,255: red,255; green,213; blue,164}, draw={rgb,255: red,121; green,57; blue,2}, shape=circle, inner sep=0pt, minimum size=0.15cm]
\tikzstyle{mid-nc}=[fill=white, draw=black, shape=circle, inner sep=0pt, minimum size=0.25cm]

\tikzstyle{di}=[draw=black, ->]
\tikzstyle{da}=[-, dashed]
\tikzstyle{da-di}=[dashed, ->]
\tikzstyle{th}=[-, thick]
\tikzstyle{th-di}=[->, thick]
\tikzstyle{th-da}=[-, dashed, thick]
\tikzstyle{gr}=[-, draw={rgb,255: red,0; green,186; blue,143}, thick]
\tikzstyle{gr-di}=[draw={rgb,255: red,0; green,186; blue,143}, ->, thick]
\tikzstyle{gr-da}=[-, dashed, draw={rgb,255: red,0; green,186; blue,143}, thick]
\tikzstyle{pu}=[-, draw={rgb,255: red,124; green,95; blue,239}, thick]
\tikzstyle{pu-da}=[-, dashed, draw={rgb,255: red,124; green,95; blue,239}, thick]
\tikzstyle{re}=[-, draw={rgb,255: red,255; green,5; blue,80}, thick]
\tikzstyle{re-di}=[->, draw={rgb,255: red,255; green,5; blue,80}, thick]
\tikzstyle{re-da}=[-, draw={rgb,255: red,255; green,5; blue,80}, dashed, thick]
\tikzstyle{bl}=[-, draw={rgb,255: red,14; green,115; blue,177}, thick]
\tikzstyle{bl-di}=[draw={rgb,255: red,14; green,115; blue,177}, ->, thick]
\tikzstyle{bl-da}=[-, draw={rgb,255: red,14; green,115; blue,177}, dashed, thick]
\tikzstyle{lgr}=[-, draw={rgb,255: red,188; green,220; blue,156}, thick]
\tikzstyle{lgr-da}=[-, draw={rgb,255: red,188; green,220; blue,156}, thick, dashed]
\tikzstyle{lbl}=[-, draw={rgb,255: red,168; green,207; blue,245}, thick]
\tikzstyle{lbl-da}=[-, draw={rgb,255: red,168; green,207; blue,245}, thick, dashed]
\tikzstyle{pi}=[-, draw={rgb,255: red,208; green,160; blue,160}, thick]
\tikzstyle{pi-da}=[-, draw={rgb,255: red,208; green,160; blue,160}, thick, dashed]
\tikzstyle{sh}=[-, draw={rgb,255: red,171; green,171; blue,171}]
\tikzstyle{sh-da}=[-, draw={rgb,255: red,171; green,171; blue,171}, dashed]
\tikzstyle{fi-gr}=[-, fill={rgb,255: red,194; green,228; blue,162}, draw={rgb,255: red,194; green,228; blue,162}]
\tikzstyle{fi-bl}=[-, fill={rgb,255: red,175; green,215; blue,255}, draw={rgb,255: red,175; green,215; blue,255}]
\tikzstyle{fi-pi}=[-, fill={rgb,255: red,246; green,188; blue,188}, draw={rgb,255: red,246; green,188; blue,188}]
\tikzstyle{fi-pu}=[-, fill={rgb,255: red,230; green,203; blue,246}, draw={rgb,255: red,230; green,203; blue,246}]
\tikzstyle{fi-ye}=[-, fill={rgb,255: red,255; green,255; blue,155}, draw={rgb,255: red,255; green,255; blue,155}]
\tikzstyle{fi-or}=[-, fill={rgb,255: red,255; green,197; blue,51}, draw={rgb,255: red,255; green,197; blue,51}]
\tikzstyle{fi-dg}=[-, fill={rgb,255: red,152; green,170; blue,139}, draw={rgb,255: red,152; green,170; blue,139}]
\tikzstyle{fi-db}=[-, fill={rgb,255: red,111; green,147; blue,183}, draw={rgb,255: red,111; green,147; blue,183}]
\tikzstyle{fi-sh}=[-, fill={rgb,255: red,222; green,222; blue,222}, draw={rgb,255: red,222; green,222; blue,222}]
\tikzstyle{fi-dsh}=[-, fill={rgb,255: red,204; green,204; blue,204}, draw={rgb,255: red,204; green,204; blue,204}]
\tikzstyle{vt}=[-, very thick]
\tikzstyle{vt-di}=[->, very thick]
\tikzstyle{vt-da}=[-, dashed, very thick]
\tikzstyle{vt-gr}=[-, draw={rgb,255: red,0; green,186; blue,143}, very thick]
\tikzstyle{vt-gr-di}=[draw={rgb,255: red,0; green,186; blue,143}, ->, very thick]
\tikzstyle{vt-gr-da}=[-, dashed, draw={rgb,255: red,0; green,186; blue,143}, very thick]
\tikzstyle{vt-pu}=[-, draw={rgb,255: red,124; green,95; blue,239}, very thick]
\tikzstyle{vt-pu-da}=[-, dashed, draw={rgb,255: red,124; green,95; blue,239}, very thick]
\tikzstyle{vt-re}=[-, draw={rgb,255: red,255; green,5; blue,80}, very thick]
\tikzstyle{vt-re-di}=[->, draw={rgb,255: red,255; green,5; blue,80}, very thick]
\tikzstyle{vt-re-da}=[-, draw={rgb,255: red,255; green,5; blue,80}, dashed, very thick]
\tikzstyle{vt-bl}=[-, draw={rgb,255: red,14; green,115; blue,177}, very thick]
\tikzstyle{vt-bl-di}=[draw={rgb,255: red,14; green,115; blue,177}, ->, very thick]
\tikzstyle{vt-bl-da}=[-, draw={rgb,255: red,14; green,115; blue,177}, dashed, very thick]
\tikzstyle{vt-bl-da-di}=[draw={rgb,255: red,14; green,115; blue,177}, ->, very thick, dashed]
\tikzstyle{vt-lgr}=[-, draw={rgb,255: red,188; green,220; blue,156}, very thick]
\tikzstyle{vt-lgr-di}=[draw={rgb,255: red,188; green,220; blue,156}, ->, very thick]
\tikzstyle{vt-lgr-da}=[-, draw={rgb,255: red,188; green,220; blue,156}, very thick, dashed]
\tikzstyle{vt-lbl}=[-, draw={rgb,255: red,168; green,207; blue,245}, very thick]
\tikzstyle{vt-lbl-di}=[->, draw={rgb,255: red,168; green,207; blue,245}, very thick]
\tikzstyle{vt-lbl-da}=[-, draw={rgb,255: red,168; green,207; blue,245}, very thick, dashed]
\tikzstyle{vt-lpp}=[-, draw={rgb,255: red,223; green,199; blue,240}, very thick]
\tikzstyle{vt-pi}=[-, draw={rgb,255: red,236; green,151; blue,151}, very thick]
\tikzstyle{vt-pi-da}=[-, draw={rgb,255: red,236; green,151; blue,151}, very thick, dashed]
\tikzstyle{vt-ye}=[-, draw={rgb,255: red,241; green,241; blue,115}, very thick]
\tikzstyle{vt-or}=[-, draw={rgb,255: red,245; green,150; blue,32}, very thick]
\tikzstyle{vt-sh}=[-, fill=none, draw={rgb,255: red,171; green,171; blue,171}, very thick]
\tikzstyle{vt-sh-di}=[draw={rgb,255: red,171; green,171; blue,171}, very thick, ->]
\tikzstyle{vt-sh-da}=[-, fill=none, draw={rgb,255: red,171; green,171; blue,171}, very thick, dashed]
\tikzstyle{fi-bx}=[-, fill={rgb,255: red,172; green,229; blue,232}, draw={rgb,255: red,14; green,115; blue,177}, thick]
\tikzstyle{ut-gr}=[-, draw={rgb,255: red,0; green,186; blue,143}, line width=3pt]
\tikzstyle{ut-gr-di}=[draw={rgb,255: red,0; green,186; blue,143}, ->, line width=3pt]
\tikzstyle{ut-ig}=[-, draw={rgb,255: red,255; green,5; blue,80}, line width=3pt]
\tikzstyle{ut-ig-di}=[draw={rgb,255: red,255; green,5; blue,80}, line width=3pt, ->]
\tikzstyle{fi-lbl}=[-, fill={rgb,255: red,213; green,237; blue,255}, draw={rgb,255: red,213; green,237; blue,255}]
\tikzstyle{fi-vlbl}=[-, fill={rgb,255: red,229; green,251; blue,255}, draw={rgb,255: red,229; green,251; blue,255}]
\tikzstyle{op-sh}=[-, draw={rgb,255: red,113; green,113; blue,113}, fill={rgb,255: red,113; green,113; blue,113}, opacity=0.5]

%% file: figures/TFT0.tikz
\begin{tikzpicture}
	\begin{pgfonlayer}{nodelayer}
		\node [style=none] (116) at (0, -4.5) {};
		\node [style=none] (117) at (0, 4.5) {};
		\node [style=none] (118) at (4.5, 0) {};
		\node [style=none] (119) at (1, 0) {};
		\node [style=none] (49) at (0, -4.5) {};
		\node [style=none] (50) at (0, 4.5) {};
		\node [style=none] (51) at (-4.5, 0) {};
		\node [style=none] (52) at (4.5, 0) {};
		\node [style=none] (95) at (1, 0) {};
		\node [style=bdot] (96) at (-1, 0) {};
		\node [style=none] (97) at (-2.9, -3.45) {};
		\node [style=none] (98) at (-2.975, -3.4) {};
		\node [style=none] (99) at (-2.975, 3.375) {};
		\node [style=none] (100) at (-2.9, 3.425) {};
		\node [style=none] (101) at (0.275, 2.9) {};
		\node [style=none] (102) at (0.25, 3) {};
		\node [style=none] (103) at (0.225, -3.075) {};
		\node [style=none] (104) at (0.25, -2.975) {};
		\node [style=bdot] (105) at (1, 0) {};
		\node [style=bdot] (106) at (-4.5, 0) {};
		\node [style=none] (107) at (5.25, 0) {$N$};
		\node [style=none] (108) at (-3.5, 4) {$M'$};
		\node [style=none] (109) at (-3.5, -4) {$M$};
		\node [style=none] (110) at (-0.5, 3.25) {$Y$};
		\node [style=none] (111) at (-0.5, -3.25) {$X$};
		\node [style=none] (112) at (0, -4.5) {};
		\node [style=none] (113) at (0, 4.5) {};
		\node [style=none] (114) at (-4.5, 0) {};
		\node [style=none] (115) at (1, 0) {};
		\node [style=none] (120) at (-5.55, 0) {$\psi^{M,M'}_{i}$};
		\node [style=none] (121) at (-1.75, 0) {$\varphi^{A}_{r,s}$};
		\node [style=none] (122) at (2, 0) {$\phi^{X,Y}_{j,k}$};
		\node [style=none] (123) at (0, -5) {$\alpha$};
		\node [style=none] (124) at (0, 5) {$\beta$};
		\node [style=dot] (125) at (0, 4.5) {};
		\node [style=dot] (126) at (0, -4.5) {};
	\end{pgfonlayer}
	\begin{pgfonlayer}{edgelayer}
		\draw [style=fi-gr] (113.center)
			 to [bend right=45] (114.center)
			 to [bend right=45] (112.center)
			 to [in=-90, out=90] (115.center)
			 to [in=-90, out=90] cycle;
		\draw [style=fi-bl] (119.center)
			 to [in=-90, out=90] (117.center)
			 to [bend left=45] (118.center)
			 to [bend left=45] (116.center)
			 to [in=-90, out=90] cycle;
		\draw [style=vt-bl, in=-90, out=90] (49.center) to (95.center);
		\draw [style=vt-bl, in=-90, out=90] (95.center) to (50.center);
		\draw [style=vt-re-di, bend left=45] (50.center) to (52.center);
		\draw [style=vt-re, bend left=45] (52.center) to (49.center);
		\draw [style=vt-re, bend left=45] (49.center) to (51.center);
		\draw [style=vt-re, bend left=45] (51.center) to (50.center);
		\draw [style=vt-re-di] (97.center) to (98.center);
		\draw [style=vt-re-di] (99.center) to (100.center);
		\draw [style=vt-bl-di] (101.center) to (102.center);
		\draw [style=vt-bl-di] (103.center) to (104.center);
	\end{pgfonlayer}
\end{tikzpicture}

%% file: figures/TFT4.tikz
\begin{tikzpicture}
	\begin{pgfonlayer}{nodelayer}
		\node [style=none] (1) at (-2.5, 0) {};
		\node [style=none] (3) at (2.5, 0) {};
		\node [style=none] (4) at (0, -0.65) {};
		\node [style=none] (5) at (0, 0.65) {};
		\node [style=none] (0) at (0, 2.5) {};
		\node [style=none] (2) at (0, -2.5) {};
		\node [style=bdot] (6) at (-1.5, -0.5) {};
		\node [style=bdot] (7) at (-0.75, 1.75) {};
		\node [style=bdot] (8) at (0.75, 1.75) {};
		\node [style=bdot] (9) at (-0.75, -1.75) {};
		\node [style=bdot] (10) at (0.75, -1.75) {};
		\node [style=none] (11) at (-1, 1.325) {\small$s$};
		\node [style=none] (12) at (1, 1.325) {\small$k$};
		\node [style=none] (13) at (-1, -1.325) {\small$r$};
		\node [style=none] (14) at (1, -1.325) {\small$j$};
		\node [style=none] (15) at (-1.75, -0.75) {\small$i$};
	\end{pgfonlayer}
	\begin{pgfonlayer}{edgelayer}
		\draw [style=vt-sh] (1.center)
			 to [in=180, out=-90, looseness=0.50] (4.center)
			 to [in=-90, out=0, looseness=0.50] (3.center);
		\draw [style=vt-sh-da] (1.center)
			 to [in=-180, out=90, looseness=0.50] (5.center)
			 to [in=90, out=0, looseness=0.50] (3.center);
		\draw [style=vt] (0.center)
			 to [bend left=45] (3.center)
			 to [bend left=45] (2.center)
			 to [bend left=45] (1.center)
			 to [bend left=45] cycle;
	\end{pgfonlayer}
\end{tikzpicture}

%% file: figures/TFT1a.tikz
\begin{tikzpicture}
	\begin{pgfonlayer}{nodelayer}
		\node [style=none] (5) at (0, 1.9) {};
		\node [style=bdot] (6) at (-3.325, -1.475) {};
		\node [style=bdot] (7) at (-2, 3.25) {};
		\node [style=bdot] (8) at (2, 3.75) {};
		\node [style=bdot] (9) at (-2, -3.75) {};
		\node [style=bdot] (10) at (2, -3.25) {};
		\node [style=none] (11) at (-2.3, 2.625) {\tiny$s$};
		\node [style=none] (12) at (2.25, 3.075) {\tiny$k$};
		\node [style=none] (13) at (-2.25, -3.075) {\tiny$r$};
		\node [style=none] (14) at (2.425, -2.65) {\tiny$j$};
		\node [style=none] (15) at (-3.575, -1.725) {\tiny$i$};
		\node [style=dot] (25) at (-1, -0.25) {};
		\node [style=none] (28) at (-3.425, 0.925) {};
		\node [style=none] (29) at (-3.35, 0.975) {};
		\node [style=dot] (35) at (-2.75, -1.175) {};
		\node [style=none] (36) at (3.5, 0.5) {\tiny$N$};
		\node [style=none] (37) at (-3.7, 0.25) {\tiny$M'$};
		\node [style=none] (38) at (0.3, -1.125) {\tiny$M$};
		\node [style=none] (39) at (-2.3, 0.45) {\tiny$Y$};
		\node [style=none] (40) at (2.5, -0.5) {\tiny$X$};
		\node [style=none] (45) at (-2.3, -0.9) {\tiny$\psi_{i}$};
		\node [style=none] (46) at (-0.325, -0.475) {\tiny$\varphi_{r,s}$};
		\node [style=none] (47) at (1.775, 0.5) {\tiny$\phi_{j,k}$};
		\node [style=none] (48) at (3.475, -1.175) {\tiny$\alpha$};
		\node [style=none] (49) at (-2.825, 1.4) {\tiny$\beta$};
		\node [style=none] (52) at (-4.25, 0) {};
		\node [style=none] (53) at (4.25, 0) {};
		\node [style=none] (54) at (0, -1.4) {};
		\node [style=none] (55) at (0, 1.4) {};
		\node [style=bdot] (56) at (1, 0.5) {};
		\node [style=dot] (58) at (1, 0.5) {};
		\node [style=dot] (63) at (3, -1.1) {};
		\node [style=dot] (64) at (-2.975, 1.1) {};
		\node [style=none] (65) at (-2.1, 0.825) {};
		\node [style=none] (66) at (-2.175, 0.85) {};
		\node [style=none] (67) at (1.95, -0.425) {};
		\node [style=none] (68) at (2, -0.5) {};
		\node [style=none] (69) at (-0.35, -1.4) {};
		\node [style=none] (70) at (-0.2, -1.4) {};
		\node [style=none] (71) at (-1.275, 0.775) {};
		\node [style=none] (72) at (-0.8, 0.775) {};
		\node [style=none] (73) at (-1.45, 1.375) {};
		\node [style=none] (74) at (-0.95, 1.4) {};
		\node [style=none] (75) at (0.85, 1.4) {};
		\node [style=none] (76) at (1.275, 1.375) {};
		\node [style=none] (81) at (0.4, -0.45) {};
		\node [style=none] (82) at (1.25, -0.25) {};
		\node [style=none] (83) at (1.5, 0.15) {};
		\node [style=none] (84) at (2, -1.3) {};
		\node [style=none] (85) at (0, 0.775) {};
		\node [style=none] (86) at (-0.775, 0.25) {};
		\node [style=none] (87) at (0.2, 0.25) {};
		\node [style=none] (88) at (-1.05, -1.175) {};
		\node [style=none] (90) at (-1.225, 0.25) {};
		\node [style=none] (91) at (-1.5, 0.25) {};
		\node [style=none] (92) at (-1.75, 0.775) {};
		\node [style=none] (93) at (-3, 0) {};
		\node [style=none] (94) at (-3.75, 0.75) {};
		\node [style=none] (95) at (-2.25, -1.275) {};
		\node [style=none] (96) at (-3.4, -0.975) {};
		\node [style=none] (97) at (0.5, -1.375) {};
		\node [style=none] (98) at (-2.475, -0.15) {};
		\node [style=none] (102) at (0.5, -0.8) {};
		\node [style=none] (103) at (0.25, 0.775) {};
		\node [style=none] (104) at (0.5, 1.075) {};
		\node [style=none] (105) at (0.3, 1.375) {};
		\node [style=none] (106) at (0.825, 1) {};
		\node [style=none] (107) at (1.25, 0.95) {};
		\node [style=none] (109) at (2.075, 0.825) {};
		\node [style=none] (110) at (2.25, 1.25) {};
		\node [style=none] (111) at (2.25, 0.5) {};
		\node [style=none] (112) at (1.675, 0.025) {};
		\node [style=none] (113) at (3, 0) {};
		\node [style=none] (114) at (4.1, 0.4) {};
		\node [style=none] (115) at (3.75, -0.775) {};
		\node [style=none] (116) at (-0.35, 1.075) {};
		\node [style=none] (117) at (-0.175, 0.8) {};
		\node [style=none] (118) at (-1.5, 1.075) {};
		\node [style=none] (119) at (-1.35, 1.075) {};
		\node [style=none] (120) at (-0.925, 1.075) {};
		\node [style=none] (121) at (-1.5, 0.8) {};
		\node [style=none] (122) at (-2.025, 1.3) {};
		\node [style=none] (127) at (-0.25, 0.25) {};
		\node [style=none] (128) at (-1.55, -0.825) {};
		\node [style=none] (129) at (-1.7, -1.325) {};
		\node [style=none] (131) at (-1.8, -0.45) {};
		\node [style=none] (132) at (-2.7, 0.05) {};
		\node [style=none] (4) at (0, -1.9) {};
		\node [style=none] (0) at (0, 5) {};
		\node [style=none] (2) at (0, -5) {};
		\node [style=none] (1) at (-5, 0) {};
		\node [style=none] (3) at (5, 0) {};
		\node [style=none] (133) at (1, -0.075) {};
		\node [style=none] (134) at (1.025, -0.5) {};
		\node [style=none] (135) at (1.125, -1.175) {};
		\node [style=none] (136) at (1.225, -1.575) {};
		\node [style=none] (137) at (-1.025, -0.9) {};
		\node [style=none] (138) at (-1.1, -1.6) {};
		\node [style=none] (139) at (-1, -0.575) {};
	\end{pgfonlayer}
	\begin{pgfonlayer}{edgelayer}
		\draw [style=vt-lgr] (81.center) to (82.center);
		\draw [style=vt-lgr] (83.center) to (82.center);
		\draw [style=vt-lgr] (82.center) to (84.center);
		\draw [style=vt-lgr] (85.center) to (81.center);
		\draw [style=vt-lgr] (86.center) to (87.center);
		\draw [style=vt-lgr] (92.center) to (91.center);
		\draw [style=vt-lgr] (93.center) to (91.center);
		\draw [style=vt-lgr] (94.center) to (93.center);
		\draw [style=vt-lgr] (96.center) to (93.center);
		\draw [style=vt-sh-da] (1.center)
			 to [in=-180, out=90, looseness=0.75] (5.center)
			 to [in=90, out=0, looseness=0.75] (3.center);
		\draw [style=vt-bl-di] (28.center) to (29.center);
		\draw [style=vt-bl-di] (65.center) to (66.center);
		\draw [style=vt-bl-di] (68.center) to (67.center);
		\draw [style=vt-bl-di] (70.center) to (69.center);
		\draw [style=vt-bl-di] (35) to (6);
		\draw [style=vt-lgr] (91.center) to (90.center);
		\draw [style=vt-lbl] (105.center) to (104.center);
		\draw [style=vt-lbl] (104.center) to (103.center);
		\draw [style=vt-lbl] (106.center) to (104.center);
		\draw [style=vt-lbl] (110.center) to (109.center);
		\draw [style=vt-lbl] (109.center) to (107.center);
		\draw [style=vt-lbl] (109.center) to (111.center);
		\draw [style=vt-lbl] (112.center) to (111.center);
		\draw [style=vt-lbl] (111.center) to (113.center);
		\draw [style=vt-lbl] (114.center) to (113.center);
		\draw [style=vt-lbl] (113.center) to (115.center);
		\draw [style=vt-lbl] (116.center) to (55.center);
		\draw [style=vt-lbl] (116.center) to (117.center);
		\draw [style=vt-lbl] (120.center) to (116.center);
		\draw [style=vt-lbl] (118.center) to (119.center);
		\draw [style=vt-lbl] (118.center) to (121.center);
		\draw [style=vt-lbl] (122.center) to (118.center);
		\draw [style=vt-bl, in=-30, out=165, looseness=1.25] (63) to (58);
		\draw [style=vt-bl, in=-30, out=180, looseness=0.75] (71.center) to (64);
		\draw [style=vt-bl, in=150, out=0] (72.center) to (58);
		\draw [style=vt-bl, in=90, out=-180, looseness=0.75] (73.center) to (52.center);
		\draw [style=vt-bl-di, in=90, out=0, looseness=0.75] (76.center) to (53.center);
		\draw [style=vt-bl, in=360, out=180] (75.center) to (55.center);
		\draw [style=vt-bl, in=180, out=0] (74.center) to (55.center);
		\draw [style=vt-bl, in=-90, out=0, looseness=0.75] (54.center) to (53.center);
		\draw [style=vt-lgr] (127.center) to (25);
		\draw [style=vt-lgr] (128.center) to (81.center);
		\draw [style=vt-lgr] (128.center) to (129.center);
		\draw [style=vt-lgr] (131.center) to (25);
		\draw [style=vt-lgr] (131.center) to (128.center);
		\draw [style=vt-lgr] (132.center) to (131.center);
		\draw [style=vt-bl, in=180, out=-90, looseness=0.75] (52.center) to (54.center);
		\draw [style=vt-sh] (1.center)
			 to [in=180, out=-90, looseness=0.75] (4.center)
			 to [in=-90, out=0, looseness=0.75] (3.center);
		\draw [style=vt] (0.center)
			 to [bend left=45] (3.center)
			 to [bend left=45] (2.center)
			 to [bend left=45] (1.center)
			 to [bend left=45] cycle;
		\draw [style=vt-bl-di, bend right=15] (25) to (7);
		\draw [style=vt-bl-di, bend left=15] (58) to (8);
		\draw [style=vt-bl, in=90, out=-90] (58) to (133.center);
		\draw [style=vt-bl] (134.center) to (135.center);
		\draw [style=vt-bl-di, in=120, out=-75, looseness=0.75] (136.center) to (10);
		\draw [style=vt-bl] (137.center) to (88.center);
		\draw [style=vt-bl] (25) to (139.center);
		\draw [style=vt-bl-di, in=60, out=-105] (138.center) to (9);
	\end{pgfonlayer}
\end{tikzpicture}

%% file: figures/TFT2.tikz
\begin{tikzpicture}
	\begin{pgfonlayer}{nodelayer}
		\node [style=none] (14) at (-1.65, 1.975) {};
		\node [style=none] (4) at (-0.5, 0.5) {};
		\node [style=none] (5) at (0.25, -0.5) {};
		\node [style=none] (6) at (2.875, -3.5) {};
		\node [style=none] (7) at (-2, 2.5) {};
		\node [style=none] (8) at (0, 0) {$\phi_{j,k}$};
		\node [style=none] (9) at (0.75, -0.5) {};
		\node [style=none] (10) at (2, -0.75) {};
		\node [style=none] (11) at (2, 3) {};
		\node [style=none] (12) at (0, -4) {};
		\node [style=none] (13) at (-2, -0.75) {};
		\node [style=none] (15) at (0.875, -1.25) {};
		\node [style=none] (16) at (-4, -2) {};
		\node [style=none] (17) at (0.25, -2.5) {};
		\node [style=none] (18) at (1.5, -4) {};
		\node [style=none] (19) at (0, -1.25) {};
		\node [style=none] (20) at (0, -2.25) {};
		\node [style=none] (21) at (0, -3) {};
		\node [style=none] (22) at (-0.25, -0.5) {};
		\node [style=none] (23) at (2.75, -1.25) {};
		\node [style=none] (24) at (1.75, -2.2) {};
		\node [style=none] (25) at (4, -2.25) {};
		\node [style=none] (26) at (3, -0.5) {};
		\node [style=none] (27) at (2.25, 0.25) {};
		\node [style=none] (28) at (4, 0.25) {};
		\node [style=none] (29) at (-0.85, 1) {};
		\node [style=none] (30) at (0.25, 1.75) {};
		\node [style=none] (31) at (1.75, 0.8) {};
		\node [style=none] (32) at (0, 3) {};
		\node [style=none] (33) at (2.275, -2.825) {};
		\node [style=none] (34) at (2.2, -2.75) {};
		\node [style=none] (0) at (-1.5, 0.5) {};
		\node [style=none] (1) at (0.5, 0.5) {};
		\node [style=none] (2) at (-0.75, -0.5) {};
		\node [style=none] (3) at (1.25, -0.5) {};
		\node [style=none] (35) at (-1.25, 2.5) {$Y$};
		\node [style=none] (36) at (3.25, -3) {$X$};
		\node [style=none] (37) at (2.5, 2.5) {$k$};
		\node [style=none] (38) at (-0.5, -3.5) {$j$};
		\node [style=none] (39) at (-2.5, 0.75) {};
		\node [style=none] (40) at (-4, 1) {};
	\end{pgfonlayer}
	\begin{pgfonlayer}{edgelayer}
		\draw [style=vt-lgr, in=-135, out=30, looseness=0.75] (39.center) to (14.center);
		\draw [style=vt-lgr, in=-150, out=30, looseness=1.25] (17.center) to (15.center);
		\draw [style=vt-bl-di] (10.center) to (11.center);
		\draw [style=vt-bl, in=-90, out=-45, looseness=2.00] (9.center) to (10.center);
		\draw [style=vt-lgr, in=15, out=-165, looseness=1.25] (13.center) to (16.center);
		\draw [style=vt-lgr, in=-150, out=0] (13.center) to (17.center);
		\draw [style=vt-lgr, in=120, out=-75] (17.center) to (18.center);
		\draw [style=vt-bl-di] (21.center) to (12.center);
		\draw [style=vt-bl] (19.center) to (20.center);
		\draw [style=vt-bl, in=90, out=-30] (22.center) to (19.center);
		\draw [style=vt-lbl, in=-165, out=15, looseness=1.25] (24.center) to (23.center);
		\draw [style=vt-lbl, in=135, out=15] (23.center) to (25.center);
		\draw [style=vt-lbl, in=90, out=-120] (26.center) to (23.center);
		\draw [style=vt-lbl, in=165, out=-60] (27.center) to (26.center);
		\draw [style=vt-lbl, in=0, out=-135] (28.center) to (26.center);
		\draw [style=vt-lbl, in=-165, out=30] (29.center) to (30.center);
		\draw [style=vt-lbl, in=15, out=135, looseness=0.75] (31.center) to (30.center);
		\draw [style=vt-lbl, in=-105, out=120, looseness=0.75] (30.center) to (32.center);
		\draw [style=vt-bl] (6.center) to (5.center);
		\draw [style=vt-bl-di] (4.center) to (7.center);
		\draw [style=vt-bl-di] (33.center) to (34.center);
		\draw [style=vt] (0.center)
			 to (2.center)
			 to (3.center)
			 to (1.center)
			 to cycle;
		\draw [style=vt-lgr, in=90, out=-60, looseness=0.75] (39.center) to (13.center);
		\draw [style=vt-lgr, in=15, out=-150, looseness=0.75] (39.center) to (40.center);
	\end{pgfonlayer}
\end{tikzpicture}

%% file: figures/TFT3.tikz
\begin{tikzpicture}
	\begin{pgfonlayer}{nodelayer}
		\node [style=none] (153) at (0, -4.5) {};
		\node [style=none] (154) at (0, 4.5) {};
		\node [style=none] (155) at (4.5, 0) {};
		\node [style=none] (156) at (1, 0.5) {};
		\node [style=none] (157) at (1.5, 0) {};
		\node [style=none] (158) at (1, -0.5) {};
		\node [style=none] (139) at (0, -4.5) {};
		\node [style=none] (140) at (0, 4.5) {};
		\node [style=none] (141) at (1, 0.5) {};
		\node [style=none] (142) at (0.5, 0) {};
		\node [style=none] (143) at (1, -0.5) {};
		\node [style=none] (144) at (-1, 0.5) {};
		\node [style=none] (145) at (-1.5, 0) {};
		\node [style=none] (146) at (-0.5, 0) {};
		\node [style=none] (147) at (-1, -0.5) {};
		\node [style=none] (148) at (-4.475, 0.5) {};
		\node [style=none] (149) at (-4, 0) {};
		\node [style=none] (150) at (-4.475, -0.5) {};
		\node [style=none] (151) at (-0.5, 0) {};
		\node [style=none] (152) at (0.5, 0) {};
		\node [style=none] (49) at (0, -4.5) {};
		\node [style=none] (50) at (0, 4.5) {};
		\node [style=none] (51) at (-4.5, 0) {};
		\node [style=none] (52) at (4.5, 0) {};
		\node [style=none] (97) at (-3.025, -3.475) {};
		\node [style=none] (98) at (-3.1, -3.425) {};
		\node [style=none] (99) at (-3.1, 3.4) {};
		\node [style=none] (100) at (-3.025, 3.45) {};
		\node [style=none] (101) at (0.35, 2.9) {};
		\node [style=none] (102) at (0.3, 3) {};
		\node [style=none] (103) at (0.3, -3.075) {};
		\node [style=none] (104) at (0.35, -2.975) {};
		\node [style=none] (107) at (5.25, 0) {$N$};
		\node [style=none] (108) at (-3.5, 4) {$M'$};
		\node [style=none] (109) at (-3.5, -4) {$M$};
		\node [style=none] (110) at (-0.5, 3.25) {$Y$};
		\node [style=none] (111) at (-0.5, -3.25) {$X$};
		\node [style=none] (123) at (0, -5) {$\alpha$};
		\node [style=none] (124) at (0, 5) {$\beta$};
		\node [style=none] (125) at (1, 0.5) {};
		\node [style=none] (126) at (0.5, 0) {};
		\node [style=none] (127) at (1.5, 0) {};
		\node [style=none] (128) at (1, -0.5) {};
		\node [style=none] (129) at (-1, 0.5) {};
		\node [style=none] (130) at (-1.5, 0) {};
		\node [style=none] (131) at (-0.5, 0) {};
		\node [style=none] (132) at (-1, -0.5) {};
		\node [style=dot] (133) at (0, 4.5) {};
		\node [style=dot] (134) at (0, -4.5) {};
		\node [style=none] (135) at (-4.475, 0.5) {};
		\node [style=none] (137) at (-4, 0) {};
		\node [style=none] (138) at (-4.475, -0.5) {};
	\end{pgfonlayer}
	\begin{pgfonlayer}{edgelayer}
		\draw [style=fi-bl] (156.center)
			 to [bend left=45] (157.center)
			 to [bend left=45] (158.center)
			 to [in=90, out=-90, looseness=0.75] (153.center)
			 to [bend right=45] (155.center)
			 to [bend right=45] (154.center)
			 to [in=90, out=-90, looseness=0.75] cycle;
		\draw [style=fi-gr] (151.center)
			 to [bend left=45] (147.center)
			 to [bend left=45] (145.center)
			 to [bend left=45] (144.center)
			 to [bend left=45] (146.center)
			 to (142.center)
			 to [bend left=45] (141.center)
			 to [in=-90, out=90, looseness=0.75] (140.center)
			 to [bend right=45] (148.center)
			 to [bend left=45] (149.center)
			 to [bend left=45] (150.center)
			 to [bend right=45] (139.center)
			 to [in=-90, out=90, looseness=0.75] (143.center)
			 to [bend left=45] (152.center)
			 to cycle;
		\draw [style=vt-bl, in=90, out=-90, looseness=0.75] (50.center) to (125.center);
		\draw [style=vt-bl, in=90, out=-90, looseness=0.75] (128.center) to (49.center);
		\draw [style=vt-re-di, bend left=45] (50.center) to (52.center);
		\draw [style=vt-re, bend left=45] (52.center) to (49.center);
		\draw [style=vt-re-di] (97.center) to (98.center);
		\draw [style=vt-re-di] (99.center) to (100.center);
		\draw [style=vt-bl-di] (101.center) to (102.center);
		\draw [style=vt-bl-di] (103.center) to (104.center);
		\draw [style=vt] (125.center)
			 to [bend left=45] (127.center)
			 to [bend left=45] (128.center)
			 to [bend left=45] (126.center)
			 to [bend left=45] cycle;
		\draw [style=vt] (129.center)
			 to [bend left=45] (131.center)
			 to [bend left=45] (132.center)
			 to [bend left=45] (130.center)
			 to [bend left=45] cycle;
		\draw [style=vt-re, bend left=45] (135.center) to (50.center);
		\draw [style=vt-re, bend left=45] (49.center) to (138.center);
		\draw [style=vt, bend left=45] (135.center) to (137.center);
		\draw [style=vt, bend left=45] (137.center) to (138.center);
	\end{pgfonlayer}
\end{tikzpicture}

%% file: figures/WS1.tikz
\begin{tikzpicture}
	\begin{pgfonlayer}{nodelayer}
		\node [style=none] (71) at (-0.25, -2.25) {};
		\node [style=none] (72) at (-3.75, -2) {};
		\node [style=none] (80) at (5.65, -2) {};
		\node [style=none] (81) at (-4.25, -2.75) {};
		\node [style=none] (82) at (0.5, -3.25) {};
		\node [style=none] (83) at (3.75, -2.5) {};
		\node [style=none] (84) at (6.125, -2.925) {};
		\node [style=none] (0) at (-0.25, -1) {};
		\node [style=none] (1) at (-3, 0.25) {};
		\node [style=none] (2) at (-0.25, 1.5) {};
		\node [style=none] (3) at (2.5, 0.25) {};
		\node [style=none] (4) at (-0.25, -1) {};
		\node [style=none] (5) at (-0.25, -2.25) {};
		\node [style=none] (6) at (-0.25, -2.25) {};
		\node [style=none] (7) at (-3.75, -2) {};
		\node [style=none] (8) at (-6.75, -0.25) {};
		\node [style=none] (9) at (-7.25, 0.5) {};
		\node [style=none] (10) at (-0.25, 2.5) {};
		\node [style=none] (11) at (5.575, 2.45) {};
		\node [style=none] (12) at (6.075, 1.95) {};
		\node [style=none] (13) at (4, 0.25) {};
		\node [style=none] (14) at (5.75, -1.25) {};
		\node [style=none] (15) at (5.65, -2) {};
		\node [style=none] (196) at (5.75, -1.25) {};
		\node [style=none] (197) at (6.25, -0.5) {};
		\node [style=none] (198) at (5.65, -2) {};
		\node [style=none] (199) at (6.125, -2.925) {};
		\node [style=none] (200) at (6.75, -1.5) {};
		\node [style=none] (16) at (-0.25, -1) {};
		\node [style=none] (17) at (-3, 0.25) {};
		\node [style=none] (18) at (-0.25, 1.5) {};
		\node [style=none] (19) at (2.5, 0.25) {};
		\node [style=none] (20) at (-0.25, -1) {};
		\node [style=none] (21) at (-0.25, -0.25) {};
		\node [style=none] (22) at (-1.25, 0.25) {};
		\node [style=none] (23) at (-0.25, 0.75) {};
		\node [style=none] (24) at (0.75, 0.25) {};
		\node [style=none] (25) at (-0.25, -0.25) {};
		\node [style=none] (35) at (-7.25, 0.5) {};
		\node [style=none] (36) at (-0.25, 2.5) {};
		\node [style=none] (37) at (-7.575, 0.975) {};
		\node [style=none] (38) at (-5.75, 1.25) {};
		\node [style=none] (39) at (-6.35, 1.5) {};
		\node [style=none] (52) at (6.075, 1.95) {};
		\node [style=none] (53) at (4, 0.25) {};
		\node [style=none] (54) at (5.75, -1.25) {};
		\node [style=none] (55) at (6.125, 1.475) {};
		\node [style=none] (56) at (5.125, 0.75) {};
		\node [style=none] (57) at (4.75, 0.25) {};
		\node [style=none] (58) at (6.25, -0.5) {};
		\node [style=none] (62) at (6.125, 1.475) {};
		\node [style=none] (63) at (5.125, 0.75) {};
		\node [style=none] (64) at (5.75, 0.825) {};
		\node [style=none] (85) at (-6.75, -0.25) {};
		\node [style=none] (86) at (5.575, 2.45) {};
		\node [style=none] (87) at (-3, 0.25) {};
		\node [style=none] (88) at (-0.25, 1.5) {};
		\node [style=none] (89) at (2.5, 0.25) {};
		\node [style=none] (90) at (-0.25, -1) {};
		\node [style=none] (91) at (-1.25, 0.25) {};
		\node [style=none] (92) at (-0.25, 0.75) {};
		\node [style=none] (93) at (0.75, 0.25) {};
		\node [style=none] (94) at (-0.25, -0.25) {};
		\node [style=none] (95) at (-7.25, 0.5) {};
		\node [style=none] (96) at (-0.25, 2.5) {};
		\node [style=none] (97) at (-7.575, 0.975) {};
		\node [style=none] (98) at (-5.75, 1.25) {};
		\node [style=none] (99) at (-6.35, 1.5) {};
		\node [style=none] (101) at (4, 0.25) {};
		\node [style=none] (102) at (5.75, -1.25) {};
		\node [style=none] (103) at (4.75, 0.25) {};
		\node [style=none] (104) at (6.25, -0.5) {};
		\node [style=none] (106) at (5.125, 0.75) {};
		\node [style=none] (107) at (5.75, 0.825) {};
		\node [style=none] (108) at (-0.25, -2.25) {};
		\node [style=dot] (109) at (-3.75, -2) {};
		\node [style=none] (110) at (5.65, -2) {};
		\node [style=none] (111) at (-4.25, -2.75) {};
		\node [style=none] (112) at (0.5, -3.25) {};
		\node [style=none] (113) at (3.75, -2.5) {};
		\node [style=none] (114) at (6.125, -2.925) {};
		\node [style=none] (115) at (-6.65, 1) {};
		\node [style=none] (116) at (-6.5, -0.25) {};
		\node [style=none] (117) at (6.75, -1.5) {};
		\node [style=dot] (118) at (-0.775, -0.975) {};
		\node [style=dot] (119) at (-0.75, -2.25) {};
		\node [style=none] (120) at (-0.925, -0.1) {};
		\node [style=none] (122) at (-1.5, 0.5) {};
		\node [style=none] (124) at (5.05, -0.175) {};
		\node [style=none] (125) at (4.65, 0.325) {};
		\node [style=none] (126) at (-4.275, 1.225) {};
		\node [style=none] (127) at (-4.25, 1.25) {};
		\node [style=none] (128) at (-0.6, -2.8) {};
		\node [style=none] (129) at (-0.65, -2.725) {};
		\node [style=none] (130) at (-5, -0.5) {};
		\node [style=none] (131) at (-4.95, -0.525) {};
		\node [style=none] (132) at (-3.825, -2.325) {};
		\node [style=none] (133) at (-3.8, -2.275) {};
		\node [style=none] (134) at (-2.25, -1.95) {};
		\node [style=none] (135) at (-2.225, -1.95) {};
		\node [style=none] (138) at (2.75, -1.825) {};
		\node [style=none] (139) at (2.725, -1.825) {};
		\node [style=none] (140) at (-4.55, -2.25) {\tiny $M_1$};
		\node [style=none] (141) at (-4.25, -0.5) {\tiny $M_2$};
		\node [style=none] (142) at (-3.75, 1) {\tiny $X_1$};
		\node [style=none] (143) at (0.425, 1.775) {\tiny $X_6$};
		\node [style=none] (144) at (3.525, 0.25) {\tiny $X_2$};
		\node [style=none] (146) at (3, -1.35) {\tiny $X_3$};
		\node [style=none] (147) at (-1.125, -0.5) {\tiny $X_4$};
		\node [style=none] (149) at (-2.5, -1.5) {\tiny $X_5$};
		\node [style=none] (150) at (-1.125, -2.5) {\tiny $\varphi_1$};
		\node [style=none] (151) at (-1.225, -1.25) {\tiny $\varphi_2$};
		\node [style=none] (152) at (-3.25, -2.25) {\tiny $\varphi_3$};
		\node [style=none] (164) at (5.575, 2.45) {};
		\node [style=none] (165) at (6.075, 1.95) {};
		\node [style=none] (169) at (6.075, 1.95) {};
		\node [style=none] (170) at (4, 0.25) {};
		\node [style=none] (172) at (6.125, 1.475) {};
		\node [style=none] (176) at (6.125, 1.475) {};
		\node [style=none] (177) at (5.125, 0.75) {};
		\node [style=none] (178) at (5.75, 0.825) {};
		\node [style=none] (179) at (0.425001, -0.1) {};
		\node [style=none] (180) at (1, 0.5) {};
		\node [style=none] (181) at (3.975, 1) {};
		\node [style=none] (182) at (4.35, 0.8) {};
		\node [style=none] (183) at (0.25, -0.75) {};
		\node [style=none] (184) at (0.4, -1.2) {};
		\node [style=none] (185) at (0.425, -2.025) {};
		\node [style=none] (186) at (-0.175, -3.2) {};
		\node [style=none] (187) at (0.375, -2.475) {};
		\node [style=none] (188) at (-6.65, 0.775) {};
		\node [style=none] (189) at (-6.65, 0.7) {};
		\node [style=none] (190) at (-6.6, 0.35) {};
		\node [style=none] (191) at (-6.525, -0.075) {};
		\node [style=none] (192) at (5.35, 1.95) {};
		\node [style=none] (193) at (5.35, 1.525) {};
		\node [style=none] (194) at (5.4, 1.2) {};
		\node [style=none] (195) at (5.575, 0.95) {};
		\node [style=none] (201) at (-3.75, -2) {};
		\node [style=none] (202) at (-4.25, -2.75) {};
		\node [style=none] (203) at (-6.5, -0.25) {};
		\node [style=none] (204) at (-6.65, 1) {};
		\node [style=none] (205) at (-5.75, 1.25) {};
		\node [style=none] (206) at (-6.35, 1.5) {};
		\node [style=none] (207) at (-6.65, 1) {};
	\end{pgfonlayer}
	\begin{pgfonlayer}{edgelayer}
		\draw [style=fi-gr] (81.center)
			 to [in=-90, out=-15, looseness=0.75] (72.center)
			 to [in=180, out=15] (71.center)
			 to [in=165, out=0, looseness=1.25] (80.center)
			 to [in=180, out=-90, looseness=0.75] (84.center)
			 to [in=0, out=150] (83.center)
			 to [in=0, out=-180] (82.center)
			 to [in=15, out=180] cycle;
		\draw [style=th, in=0, out=150] (114.center) to (113.center);
		\draw [style=th, in=0, out=-180] (113.center) to (112.center);
		\draw [style=fi-bl] (3.center)
			 to [in=0, out=-90, looseness=0.75] (4.center)
			 to (6.center)
			 to [in=165, out=0, looseness=1.25] (15.center)
			 to [in=255, out=90] (14.center)
			 to [in=-90, out=165] (13.center)
			 to [in=-165, out=90] (12.center)
			 to [in=15, out=105] (11.center)
			 to [in=0, out=-165] (10.center)
			 to [in=0, out=-180] (9.center)
			 to (8.center)
			 to [in=90, out=0] (7.center)
			 to [in=180, out=15] (5.center)
			 to (0.center)
			 to [in=-90, out=180, looseness=0.75] (1.center)
			 to [in=-180, out=90, looseness=0.75] (2.center)
			 to [in=90, out=0, looseness=0.75] cycle;
		\draw [style=fi-dg] (196.center)
			 to [in=-150, out=75, looseness=0.75] (197.center)
			 to [in=90, out=30] (200.center)
			 to [in=0, out=-90, looseness=0.75] (199.center)
			 to [in=-90, out=180, looseness=0.75] (198.center)
			 to [in=255, out=90] cycle;
		\draw [style=fi-gr] (22.center)
			 to [in=-180, out=90, looseness=0.75] (23.center)
			 to [in=90, out=0, looseness=0.75] (24.center)
			 to [in=0, out=-90, looseness=0.75] (25.center)
			 to (20.center)
			 to [in=-90, out=0, looseness=0.75] (19.center)
			 to [in=0, out=90, looseness=0.75] (18.center)
			 to [in=90, out=-180, looseness=0.75] (17.center)
			 to [in=180, out=-90, looseness=0.75] (16.center)
			 to (21.center)
			 to [in=-90, out=180, looseness=0.75] cycle;
		\draw [style=fi-gr] (38.center)
			 to [in=15, out=30] (39.center)
			 to [in=180, out=0] (36.center)
			 to [in=0, out=-180] (35.center)
			 to (37.center)
			 to [in=-150, out=0, looseness=0.75] cycle;
		\draw [style=fi-gr] (56.center)
			 to [in=90, out=-165] (57.center)
			 to [in=180, out=-90, looseness=0.75] (58.center)
			 to [in=75, out=-150, looseness=0.75] (54.center)
			 to [in=-90, out=165] (53.center)
			 to [in=-165, out=90] (52.center)
			 to [in=90, out=-75, looseness=0.75] (55.center)
			 to [in=30, out=-180] cycle;
		\draw [style=fi-bl] (62.center)
			 to [in=30, out=-180] (63.center)
			 to [in=180, out=15] (64.center)
			 to [in=-90, out=0] cycle;
		\draw [style=th, in=-180, out=90, looseness=0.75] (91.center) to (92.center);
		\draw [style=th, in=90, out=0, looseness=0.75] (92.center) to (93.center);
		\draw [style=bl, in=-90, out=0, looseness=0.75] (90.center) to (89.center);
		\draw [style=bl, in=0, out=90, looseness=0.75] (89.center) to (88.center);
		\draw [style=bl-di, in=-180, out=90, looseness=0.75] (87.center) to (88.center);
		\draw [style=bl, in=-180, out=0] (95.center) to (96.center);
		\draw [style=th] (95.center) to (97.center);
		\draw [style=th, in=180, out=-90, looseness=0.75] (103.center) to (104.center);
		\draw [style=th, in=75, out=-150, looseness=0.75] (104.center) to (102.center);
		\draw [style=bl, in=-90, out=165] (102.center) to (101.center);
		\draw [style=bl, in=180, out=15] (109) to (108.center);
		\draw [style=bl, in=165, out=0, looseness=1.25] (108.center) to (110.center);
		\draw [style=th, in=180, out=-90, looseness=0.75] (110.center) to (114.center);
		\draw [style=th, in=15, out=180] (112.center) to (111.center);
		\draw [style=th, in=90, out=30] (104.center) to (117.center);
		\draw [style=th, in=0, out=-90, looseness=0.75] (117.center) to (114.center);
		\draw [style=th, in=255, out=90] (110.center) to (102.center);
		\draw [style=th] (95.center) to (85.center);
		\draw [style=bl, in=-90, out=180, looseness=0.75] (90.center) to (87.center);
		\draw [style=bl, in=-180, out=90, looseness=0.75] (118) to (94.center);
		\draw [style=th, in=165, out=-75, looseness=0.75] (122.center) to (120.center);
		\draw [style=th, in=150, out=-75] (125.center) to (124.center);
		\draw [style=bl-di, in=-150, out=30, looseness=0.75] (126.center) to (127.center);
		\draw [style=bl-di] (128.center) to (129.center);
		\draw [style=bl-di] (135.center) to (134.center);
		\draw [style=bl-di] (138.center) to (139.center);
		\draw [style=th, in=15, out=105] (165.center) to (164.center);
		\draw [style=bl-di, in=90, out=-165] (169.center) to (170.center);
		\draw [style=th, in=90, out=-75, looseness=0.75] (169.center) to (172.center);
		\draw [style=bl-di, in=-180, out=30] (177.center) to (176.center);
		\draw [style=th, in=-90, out=0] (178.center) to (176.center);
		\draw [style=th, in=15, out=-105, looseness=0.75] (180.center) to (179.center);
		\draw [style=bl-da, in=-150, out=-30] (182.center) to (177.center);
		\draw [style=bl-da, in=0, out=150] (181.center) to (96.center);
		\draw [style=bl-da, in=105, out=0, looseness=0.75] (94.center) to (183.center);
		\draw [style=bl-da, in=75, out=-75, looseness=0.75] (184.center) to (185.center);
		\draw [style=bl-da, in=0, out=-105, looseness=0.75] (187.center) to (186.center);
		\draw [style=bl, in=165, out=-90, looseness=0.75] (119) to (186.center);
		\draw [style=re-da] (188.center) to (189.center);
		\draw [style=re-da, in=105, out=-90, looseness=0.50] (190.center) to (191.center);
		\draw [style=th-da, in=90, out=-165, looseness=0.75] (164.center) to (192.center);
		\draw [style=th-da, in=90, out=-90] (193.center) to (194.center);
		\draw [style=th-da, in=-180, out=-75, looseness=0.50] (195.center) to (178.center);
		\draw [style=fi-dg] (202.center)
			 to [in=-90, out=-15, looseness=0.75] (201.center)
			 to [in=0, out=90] (203.center)
			 to [in=-180, out=-75, looseness=0.75] cycle;
		\draw [style=re, in=-180, out=-75, looseness=0.75] (116.center) to (111.center);
		\draw [style=re, in=-90, out=-15, looseness=0.75] (111.center) to (109);
		\draw [style=re-di] (131.center) to (130.center);
		\draw [style=re, in=0, out=90] (109) to (85.center);
		\draw [style=re-di] (132.center) to (133.center);
		\draw [style=th, in=-180, out=0] (99.center) to (96.center);
		\draw [style=th, in=0, out=-165] (86.center) to (96.center);
		\draw [style=fi-dg] (204.center)
			 to [in=-150, out=0, looseness=0.75] (205.center)
			 to [in=15, out=30] (206.center)
			 to [in=90, out=-165] (207.center);
		\draw [style=re, in=15, out=30] (98.center) to (99.center);
		\draw [style=re-di, in=-150, out=0, looseness=0.75] (97.center) to (98.center);
		\draw [style=re, in=90, out=-165] (99.center) to (115.center);
		\draw [style=th, in=0, out=-90, looseness=0.75] (93.center) to (94.center);
		\draw [style=th, in=-90, out=180, looseness=0.75] (94.center) to (91.center);
		\draw [style=th, in=180, out=15] (106.center) to (107.center);
		\draw [style=th, in=90, out=-165] (106.center) to (103.center);
	\end{pgfonlayer}
\end{tikzpicture}

%% file: figures/WS3.tikz
\begin{tikzpicture}
	\begin{pgfonlayer}{nodelayer}
		\node [style=none] (252) at (5.25, 1.25) {};
		\node [style=none] (224) at (-1.925, 2.25) {};
		\node [style=none] (247) at (4.1, -0.25) {};
		\node [style=none] (204) at (-6.65, 1) {};
		\node [style=none] (197) at (6.25, -0.5) {};
		\node [style=none] (199) at (6.125, -2.925) {};
		\node [style=none] (35) at (-7.1, 0.5) {};
		\node [style=none] (37) at (-7.475, 0.975) {};
		\node [style=none] (39) at (-6.35, 1.5) {};
		\node [style=none] (85) at (-6.65, -0.25) {};
		\node [style=none] (86) at (5.575, 2.45) {};
		\node [style=none] (87) at (-3, 0.25) {};
		\node [style=none] (88) at (-0.25, 1.5) {};
		\node [style=none] (89) at (2.5, 0.25) {};
		\node [style=none] (90) at (-0.25, -1) {};
		\node [style=none] (91) at (-1.25, 0.25) {};
		\node [style=none] (92) at (-0.25, 0.75) {};
		\node [style=none] (93) at (0.75, 0.25) {};
		\node [style=none] (94) at (-0.25, -0.25) {};
		\node [style=none] (95) at (-7.1, 0.5) {};
		\node [style=none] (96) at (-0.25, 2.5) {};
		\node [style=none] (97) at (-7.475, 0.975) {};
		\node [style=none] (98) at (-5.75, 1.25) {};
		\node [style=none] (99) at (-6.35, 1.5) {};
		\node [style=none] (101) at (4, 0.25) {};
		\node [style=none] (102) at (5.75, -1.25) {};
		\node [style=none] (103) at (4.75, 0.25) {};
		\node [style=none] (104) at (6.25, -0.5) {};
		\node [style=none] (106) at (5.125, 0.75) {};
		\node [style=none] (107) at (5.75, 0.825) {};
		\node [style=none] (108) at (-0.25, -2.25) {};
		\node [style=dot] (109) at (-3.75, -2) {};
		\node [style=none] (110) at (5.65, -2) {};
		\node [style=none] (111) at (-4.25, -2.75) {};
		\node [style=none] (112) at (0.5, -3.25) {};
		\node [style=none] (113) at (3.75, -2.5) {};
		\node [style=none] (114) at (6.125, -2.925) {};
		\node [style=none] (115) at (-6.65, 1) {};
		\node [style=none] (116) at (-6.5, -0.25) {};
		\node [style=none] (117) at (6.75, -1.5) {};
		\node [style=dot] (118) at (-0.775, -0.975) {};
		\node [style=dot] (119) at (-0.75, -2.25) {};
		\node [style=none] (120) at (-0.925, -0.1) {};
		\node [style=none] (122) at (-1.5, 0.5) {};
		\node [style=none] (124) at (5.05, -0.175) {};
		\node [style=none] (125) at (4.65, 0.325) {};
		\node [style=none] (126) at (-4.275, 1.225) {};
		\node [style=none] (127) at (-4.25, 1.25) {};
		\node [style=none] (128) at (-0.6, -2.8) {};
		\node [style=none] (129) at (-0.65, -2.725) {};
		\node [style=none] (130) at (-5, -0.5) {};
		\node [style=none] (131) at (-4.95, -0.525) {};
		\node [style=none] (133) at (-3.8, -2.275) {};
		\node [style=none] (134) at (-2.25, -1.95) {};
		\node [style=none] (135) at (-2.225, -1.95) {};
		\node [style=none] (138) at (2.75, -1.825) {};
		\node [style=none] (139) at (2.725, -1.825) {};
		\node [style=none] (140) at (-4.55, -2.25) {\tiny $M_1$};
		\node [style=none] (142) at (-3.75, 1) {\tiny $X_1$};
		\node [style=none] (143) at (0.425, 1.775) {\tiny $X_6$};
		\node [style=none] (144) at (3.525, 0) {\tiny $X_2$};
		\node [style=none] (146) at (3, -1.35) {\tiny $X_3$};
		\node [style=none] (147) at (-1.125, -0.5) {\tiny $X_4$};
		\node [style=none] (149) at (-2.5, -1.5) {\tiny $X_5$};
		\node [style=none] (150) at (-1.125, -2.5) {\tiny $\mathring{\varphi_1}$};
		\node [style=none] (151) at (-1.225, -1.25) {\tiny $\mathring{\varphi_2}$};
		\node [style=none] (152) at (-3.25, -2.25) {\tiny $\mathring{\varphi_3}$};
		\node [style=none] (164) at (5.575, 2.45) {};
		\node [style=none] (165) at (6.075, 1.95) {};
		\node [style=none] (169) at (6.075, 1.95) {};
		\node [style=none] (170) at (4, 0.25) {};
		\node [style=none] (172) at (6.125, 1.475) {};
		\node [style=none] (176) at (6.125, 1.475) {};
		\node [style=none] (177) at (5.125, 0.75) {};
		\node [style=none] (178) at (5.75, 0.825) {};
		\node [style=none] (181) at (3.975, 1) {};
		\node [style=none] (182) at (4.35, 0.8) {};
		\node [style=none] (183) at (0.25, -0.75) {};
		\node [style=none] (184) at (0.4, -1.2) {};
		\node [style=none] (185) at (0.425, -2.025) {};
		\node [style=none] (186) at (-0.175, -3.2) {};
		\node [style=none] (187) at (0.025, -3.15) {};
		\node [style=none] (188) at (-6.65, 0.775) {};
		\node [style=none] (189) at (-6.65, 0.7) {};
		\node [style=none] (190) at (-6.6, 0.35) {};
		\node [style=none] (191) at (-6.525, -0.075) {};
		\node [style=none] (192) at (5.35, 1.95) {};
		\node [style=none] (193) at (5.35, 1.525) {};
		\node [style=none] (194) at (5.4, 1.2) {};
		\node [style=none] (195) at (5.575, 0.95) {};
		\node [style=none] (210) at (-6.2, -2.85) {};
		\node [style=none] (213) at (-8.025, 1.3) {};
		\node [style=none] (214) at (-5.25, 0.25) {};
		\node [style=none] (215) at (-6, 0.6) {};
		\node [style=none] (216) at (-5, 0.875) {};
		\node [style=none] (217) at (-5.6, -0.325) {};
		\node [style=none] (218) at (-3.5, 0) {};
		\node [style=none] (219) at (-3.875, -1.25) {};
		\node [style=none] (220) at (-4.75, 0.175) {};
		\node [style=none] (132) at (-3.825, -2.325) {};
		\node [style=none] (225) at (-2.25, 1.1) {};
		\node [style=none] (226) at (-1.65, 1.75) {};
		\node [style=none] (227) at (-4.25, -0.5) {\tiny $M_2$};
		\node [style=none] (228) at (1.5, 1.2) {};
		\node [style=none] (229) at (2.25, 1.5) {};
		\node [style=none] (230) at (3, 2.225) {};
		\node [style=none] (231) at (3.3, 1.4) {};
		\node [style=none] (232) at (4, 0.5) {};
		\node [style=none] (233) at (5.05, 1.625) {};
		\node [style=none] (234) at (5.75, 1.35) {};
		\node [style=none] (235) at (5.45, 0.8) {};
		\node [style=none] (236) at (4.25, 2.225) {};
		\node [style=none] (237) at (4.25, 1.325) {};
		\node [style=none] (238) at (4.525, 1) {};
		\node [style=none] (239) at (2.3, 1.95) {};
		\node [style=none] (240) at (2.65, 1.725) {};
		\node [style=none] (241) at (-1.75, -0.8) {};
		\node [style=none] (242) at (-1.75, -2.025) {};
		\node [style=none] (243) at (-1.25, -1.75) {};
		\node [style=none] (244) at (2.5, -1) {};
		\node [style=none] (245) at (1.5, -0.7) {};
		\node [style=none] (246) at (1.75, -2.05) {};
		\node [style=none] (248) at (4.75, -0.9) {};
		\node [style=none] (249) at (4.5, -1.75) {};
		\node [style=none] (250) at (4.825, 1.5) {};
		\node [style=none] (251) at (5.375, 0.95) {};
		\node [style=none] (253) at (5.525, 1.775) {};
		\node [style=none] (254) at (5.05, 1.375) {};
		\node [style=none] (255) at (4.5, 0.5) {};
		\node [style=none] (256) at (4.075, 0.725) {};
		\node [style=none] (257) at (4.525, -0.25) {};
		\node [style=none] (258) at (4.275, -0.5) {};
		\node [style=none] (259) at (5.175, -0.625) {};
		\node [style=none] (260) at (5.5, -0.4) {};
		\node [style=none] (261) at (5, -1) {};
		\node [style=none] (262) at (5.525, -0.95) {};
		\node [style=none] (263) at (5.425, -1.375) {};
		\node [style=none] (264) at (5.325, -1.7) {};
		\node [style=none] (265) at (5.175, -2.1) {};
		\node [style=none] (266) at (4.75, -2.5) {};
		\node [style=none] (267) at (4, -1.725) {};
		\node [style=none] (268) at (4.15, -2.075) {};
		\node [style=none] (269) at (2.9, -2.325) {};
		\node [style=none] (270) at (2.25, -1.95) {};
		\node [style=none] (271) at (2.5, -2.725) {};
		\node [style=none] (272) at (1.25, -2.15) {};
		\node [style=none] (273) at (0.75, -2.875) {};
		\node [style=none] (274) at (1.3, -3.15) {};
		\node [style=none] (275) at (-0.45, -3) {};
		\node [style=none] (276) at (0.225, -2.8) {};
		\node [style=none] (277) at (0.35, -2.425) {};
		\node [style=none] (278) at (-1.5, -2.9) {};
		\node [style=none] (279) at (-2.75, -1.875) {};
		\node [style=none] (280) at (-2.5, -2.5) {};
		\node [style=none] (281) at (-3.25, -2.6) {};
		\node [style=none] (282) at (-2.25, 0.5) {};
		\node [style=none] (283) at (-2, 0) {};
		\node [style=none] (284) at (-2.25, -0.6) {};
		\node [style=none] (285) at (-2.775, 0.75) {};
		\node [style=none] (286) at (-1.75, 1.275) {};
		\node [style=none] (287) at (-1.5, -0.825) {};
		\node [style=none] (288) at (-0.675, -0.575) {};
		\node [style=none] (289) at (-0.25, -0.75) {};
		\node [style=none] (290) at (0.75, -0.9) {};
		\node [style=none] (292) at (-1, 1.425) {};
		\node [style=none] (293) at (-0.75, 1) {};
		\node [style=none] (294) at (-1, 0.575) {};
		\node [style=none] (295) at (0.25, 1) {};
		\node [style=none] (296) at (0.275, 0.675) {};
		\node [style=none] (297) at (1, 1.35) {};
		\node [style=none] (298) at (1.65, 0.425) {};
		\node [style=none] (299) at (0.975, 0.975) {};
		\node [style=none] (300) at (0.375, -0.825) {};
		\node [style=none] (303) at (2.475, 0.35) {};
		\node [style=none] (304) at (1.45, -0.25) {};
		\node [style=none] (305) at (2.1, -0.4) {};
		\node [style=none] (306) at (-5.75, 1) {};
		\node [style=none] (307) at (-6.25, 1.05) {};
		\node [style=none] (308) at (-5.75, 0.675) {};
		\node [style=none] (309) at (-5.275, 1.275) {};
		\node [style=none] (310) at (-5, 1.6) {};
		\node [style=none] (311) at (-4.65, 1.025) {};
		\node [style=none] (312) at (-3.5, 1.75) {};
		\node [style=none] (313) at (-3.75, 1.475) {};
		\node [style=none] (314) at (-4, 1.8) {};
		\node [style=none] (315) at (-3, 2.075) {};
		\node [style=none] (316) at (0.425, -0.1) {};
		\node [style=none] (317) at (1, 0.5) {};
	\end{pgfonlayer}
	\begin{pgfonlayer}{edgelayer}
		\draw [style=lgr] (292.center) to (293.center);
		\draw [style=lgr] (293.center) to (294.center);
		\draw [style=lgr] (293.center) to (295.center);
		\draw [style=lgr] (295.center) to (296.center);
		\draw [style=lgr] (282.center) to (283.center);
		\draw [style=lgr] (284.center) to (283.center);
		\draw [style=lgr] (285.center) to (282.center);
		\draw [style=lgr] (282.center) to (286.center);
		\draw [style=lgr] (283.center) to (120.center);
		\draw [style=lgr] (278.center) to (129.center);
		\draw [style=lgr] (272.center) to (273.center);
		\draw [style=lgr] (273.center) to (274.center);
		\draw [style=lgr] (275.center) to (273.center);
		\draw [style=lgr] (268.center) to (269.center);
		\draw [style=lgr] (270.center) to (269.center);
		\draw [style=lgr] (269.center) to (271.center);
		\draw [style=lgr] (258.center) to (257.center);
		\draw [style=lgr] (257.center) to (255.center);
		\draw [style=lgr] (255.center) to (254.center);
		\draw [style=lgr] (256.center) to (255.center);
		\draw [style=lgr] (250.center) to (252.center);
		\draw [style=lgr] (252.center) to (251.center);
		\draw [style=lgr] (253.center) to (252.center);
		\draw [style=lbl] (248.center) to (249.center);
		\draw [style=lbl] (226.center) to (88.center);
		\draw [style=lbl] (225.center) to (226.center);
		\draw [style=lbl] (226.center) to (224.center);
		\draw [style=lbl] (241.center) to (243.center);
		\draw [style=lbl] (243.center) to (242.center);
		\draw [style=lbl] (243.center) to (108.center);
		\draw [style=lbl] (244.center) to (246.center);
		\draw [style=lbl] (244.center) to (245.center);
		\draw [style=lbl] (244.center) to (247.center);
		\draw [style=lbl] (215.center) to (214.center);
		\draw [style=lbl] (214.center) to (217.center);
		\draw [style=lbl] (214.center) to (218.center);
		\draw [style=lbl] (218.center) to (87.center);
		\draw [style=lbl] (216.center) to (220.center);
		\draw [style=lbl] (218.center) to (219.center);
		\draw [style=lbl] (228.center) to (229.center);
		\draw [style=lbl, in=165, out=60, looseness=0.75] (229.center) to (230.center);
		\draw [style=lbl-da, in=90, out=-15, looseness=0.75] (230.center) to (231.center);
		\draw [style=lbl] (229.center) to (232.center);
		\draw [style=lbl-da, in=-180, out=120] (237.center) to (236.center);
		\draw [style=lbl-da, in=-165, out=-45, looseness=0.75] (238.center) to (235.center);
		\draw [style=lbl, in=15, out=-120] (234.center) to (235.center);
		\draw [style=fi-gr] (35.center) to (37.center);
		\draw [style=th, in=-180, out=90, looseness=0.75] (91.center) to (92.center);
		\draw [style=th, in=90, out=0, looseness=0.75] (92.center) to (93.center);
		\draw [style=bl, in=15, out=30] (98.center) to (99.center);
		\draw [style=th, in=75, out=-150, looseness=0.75] (104.center) to (102.center);
		\draw [style=th, in=180, out=-90, looseness=0.75] (110.center) to (114.center);
		\draw [style=bl-da, in=90, out=-165] (99.center) to (115.center);
		\draw [style=bl-da, in=-180, out=-75, looseness=0.75] (116.center) to (111.center);
		\draw [style=th, in=90, out=30] (104.center) to (117.center);
		\draw [style=th, in=0, out=-90, looseness=0.75] (117.center) to (114.center);
		\draw [style=th, in=255, out=90] (110.center) to (102.center);
		\draw [style=th, in=165, out=-75, looseness=0.75] (122.center) to (120.center);
		\draw [style=th, in=150, out=-75] (125.center) to (124.center);
		\draw [style=bl-di, in=-150, out=30, looseness=0.75] (126.center) to (127.center);
		\draw [style=bl-di] (128.center) to (129.center);
		\draw [style=bl-di] (131.center) to (130.center);
		\draw [style=bl-di] (132.center) to (133.center);
		\draw [style=bl-di] (135.center) to (134.center);
		\draw [style=bl-di] (138.center) to (139.center);
		\draw [style=th, in=15, out=105] (165.center) to (164.center);
		\draw [style=th, in=90, out=-75, looseness=0.75] (169.center) to (172.center);
		\draw [style=th, in=-90, out=0] (178.center) to (176.center);
		\draw [style=bl-da, in=-150, out=-30] (182.center) to (177.center);
		\draw [style=bl-da, in=105, out=0, looseness=0.75] (94.center) to (183.center);
		\draw [style=bl-da, in=75, out=-75, looseness=0.75] (184.center) to (185.center);
		\draw [style=bl-da, in=0, out=-105, looseness=0.75] (187.center) to (186.center);
		\draw [style=bl, in=165, out=-90, looseness=0.75] (119) to (186.center);
		\draw [style=bl-da] (188.center) to (189.center);
		\draw [style=bl-da, in=105, out=-90, looseness=0.50] (190.center) to (191.center);
		\draw [style=th-da, in=90, out=-165, looseness=0.75] (164.center) to (192.center);
		\draw [style=th-da, in=90, out=-90] (193.center) to (194.center);
		\draw [style=th-da, in=-180, out=-75, looseness=0.50] (195.center) to (178.center);
		\draw [style=th, in=-165, out=15] (210.center) to (111.center);
		\draw [style=th, in=180, out=0] (213.center) to (99.center);
		\draw [style=th] (210.center)
			 to [in=180, out=-135, looseness=0.75] (213.center)
			 to [in=45, out=0, looseness=0.50] cycle;
		\draw [style=bl, in=-90, out=-15, looseness=0.75] (111.center) to (109);
		\draw [style=bl, in=0, out=90] (109) to (85.center);
		\draw [style=bl-di, in=90, out=-165] (169.center) to (170.center);
		\draw [style=bl-di, in=-180, out=30] (177.center) to (176.center);
		\draw [style=bl-di, in=-180, out=90, looseness=0.75] (87.center) to (88.center);
		\draw [style=th, in=180, out=15] (106.center) to (107.center);
		\draw [style=th, in=90, out=-165] (106.center) to (103.center);
		\draw [style=lbl, in=0, out=120] (233.center) to (236.center);
		\draw [style=bl-da, in=0, out=150, looseness=0.75] (239.center) to (96.center);
		\draw [style=bl-da] (240.center) to (181.center);
		\draw [style=th, in=0, out=-165] (86.center) to (96.center);
		\draw [style=lgr] (257.center) to (259.center);
		\draw [style=lgr, in=150, out=45, looseness=0.75] (259.center) to (260.center);
		\draw [style=lgr] (259.center) to (261.center);
		\draw [style=lgr-da, in=75, out=-30, looseness=0.75] (260.center) to (262.center);
		\draw [style=th, in=180, out=-90, looseness=0.75] (103.center) to (104.center);
		\draw [style=bl, in=-90, out=165] (102.center) to (101.center);
		\draw [style=lgr-da] (263.center) to (264.center);
		\draw [style=lgr-da, in=0, out=-105] (265.center) to (266.center);
		\draw [style=lgr, in=180, out=-75, looseness=0.75] (267.center) to (266.center);
		\draw [style=bl, in=165, out=0, looseness=1.25] (108.center) to (110.center);
		\draw [style=th, in=0, out=-180] (113.center) to (112.center);
		\draw [style=th, in=0, out=150] (114.center) to (113.center);
		\draw [style=bl-da] (277.center) to (276.center);
		\draw [style=lgr] (280.center) to (279.center);
		\draw [style=lgr] (280.center) to (278.center);
		\draw [style=lgr] (281.center) to (280.center);
		\draw [style=bl, in=180, out=15] (109) to (108.center);
		\draw [style=th, in=15, out=180] (112.center) to (111.center);
		\draw [style=lgr] (288.center) to (289.center);
		\draw [style=lgr] (289.center) to (90.center);
		\draw [style=lgr] (289.center) to (290.center);
		\draw [style=bl, in=-180, out=90, looseness=0.75] (118) to (94.center);
		\draw [style=th, in=-90, out=180, looseness=0.75] (94.center) to (91.center);
		\draw [style=th, in=0, out=-90, looseness=0.75] (93.center) to (94.center);
		\draw [style=lgr] (299.center) to (298.center);
		\draw [style=lgr] (298.center) to (303.center);
		\draw [style=lgr] (304.center) to (305.center);
		\draw [style=lgr] (304.center) to (298.center);
		\draw [style=lgr] (304.center) to (300.center);
		\draw [style=lgr] (295.center) to (299.center);
		\draw [style=lgr] (299.center) to (297.center);
		\draw [style=bl, in=0, out=90, looseness=0.75] (89.center) to (88.center);
		\draw [style=bl, in=-90, out=180, looseness=0.75] (90.center) to (87.center);
		\draw [style=bl, in=-90, out=0, looseness=0.75] (90.center) to (89.center);
		\draw [style=lgr] (306.center) to (308.center);
		\draw [style=lgr] (306.center) to (307.center);
		\draw [style=lgr] (309.center) to (306.center);
		\draw [style=lgr] (309.center) to (310.center);
		\draw [style=lgr] (309.center) to (311.center);
		\draw [style=lgr] (314.center) to (312.center);
		\draw [style=lgr] (312.center) to (313.center);
		\draw [style=lgr] (312.center) to (315.center);
		\draw [style=bl, in=-150, out=0, looseness=0.75] (97.center) to (98.center);
		\draw [style=bl, in=-180, out=0] (95.center) to (96.center);
		\draw [style=th, in=-180, out=0] (99.center) to (96.center);
		\draw [style=th, in=15, out=-105, looseness=0.75] (317.center) to (316.center);
	\end{pgfonlayer}
\end{tikzpicture}